\documentclass[12pt,oneside]{algorithms_and_scpp}
\usepackage{graphicx}
\usepackage{amsmath}

 \newtheorem{thm}{Theorem}[subsection]
 \newtheorem{cor}[thm]{Corollary}
 \newtheorem{lem}[thm]{Lemma}
 \newtheorem{prop}[thm]{Proposition}
 \theoremstyle{definition}
 \newtheorem{defn}[thm]{Definition}
 
 \theoremstyle{remark}
 \newtheorem{rem}[thm]{Remark}
 \newtheorem{ex}[thm]{Example}

 \newcommand{\Real}{\mathbb{R}}
 \newcommand{\Complex}{\mathbb{C}}

\begin{document}

\title[ \normalfont ORE REVISITED: AN ALGORITHMIC INVESTIGATION OF THE SIMPLE COMMUTATOR PROMISE PROBLEM ]
{ ORE REVISITED: AN ALGORITHMIC INVESTIGATION OF THE SIMPLE COMMUTATOR PROMISE PROBLEM }

\author{
{ \normalsize  \normalfont
by \\
JAMES L. ULRICH \\
\\
\\
\\
\\
\\
\\
\\
\\
\\
\\
\\
\\
\\
\\
\\
\\
\\
\\
\\
\\
\\
A dissertation submitted to the Graduate Faculty in Mathematics in partial fulfillment of the requirements for
the degree of Doctor of Philosophy, The City University of New York \\
%\\
%\\
2006 } \\
%\\
%November 29, 2005 -- draft version 6.   Alterations to algorithms of section 3.. Addition of new subsection to section $5$, and %additional questions for further research. Reworked discussion of complexity, after observations given by Professor Kossak. %Corrected various name mis-spellings, as brought to my attention by Professor O'Sullivan. \\
%\\
%October 25, 2005 -- draft version 5. Addition of subsection on word reduction in $S_n$. \\
%\\
%October 21, 2005 -- draft version 4. Complete example added in draft v3. \\
%\\
%October 21, 2005 -- draft version 3.  Corrections and addition of example to subsection 3.5.  Amended
%proof of proposition 5.2 to include assumption of efficient means for word reduction in $S_n$.\\
%\\
%October 17, 2005 -- draft version 2.  Corrections and clarifications to section 2. \\
%\\
%October 15, 2005 -- draft version 1. \\
%\\
}

%\pagenumbering{}

\maketitle
%%% ----------------------------------------------------------------------

\pagebreak

\pagestyle{myheadings}
\pagenumbering{roman}
\markright{}

\setcounter{page}{2}

 %\vspace{8in}
%\thanks{This work is being conducted under the advisement of Professor Michael Anshel.} \\
\vspace*{7in}
\begin{center} \copyright \, \, 2006 \\
  JAMES L. ULRICH \\
  All Rights Reserved
\end{center}  
%\author{James Ulrich}

%\address{Department of Mathematics, the Graduate Center of the City University of New
%York, New York, New York, 10016}

%\email{julrich@gc.cuny.edu}

\pagebreak

%%% ----------------------------------------------------------------------

\begin{center}
This manuscript has been read and accepted for the \\ 
Graduate Faculty in Mathematics in satisfaction of the \\
 dissertation requirement for the degree of Doctor of Philosophy. 
\end{center}

\begin{displaymath}
\begin{array}{lcl}

&   & \text{Professor Michael Anshel} \\
\line(1,0){100}   & &  \line(1,0){220}  \\
\text{Date} & &  \text{Chair of Examining Committee}  \\
& &  \\
 &  &   \\ 
 & &  \text{Professor J\'{o}zef Dodziuk} \\
\line(1,0){100}  & & \line(1,0){220}  \\
\text{Date} & & \text{Executive Officer}  \\
\end{array}
\end{displaymath}

\begin{displaymath}
\begin{array}{l}	
 \text{Professor Roman Kossak}   \\
    \line(1,0){200}   \\
	  \\
    \text{Professor Cormac O'Sullivan}   \\
    \line(1,0){200}   \\ 
      \\
	   \text{Professor Burton Randol}  \\
  \line(1,0){200} \\
     \\
   \text{Professor Vladimir Shpilrain}  \\
    \line(1,0){200}  \\
\end{array}
\end{displaymath}
\vspace{-0.6in}
\begin{center}
\hspace{-1.1in}Supervision Committee
\end{center}

\begin{center}
THE CITY UNIVERSITY OF NEW YORK
 \end{center}

%These are the acknowledgements.

\pagebreak

\begin{abstract}
\begin{center}
ORE REVISITED: AN ALGORITHMIC INVESTIGATION OF THE SIMPLE COMMUTATOR PROMISE PROBLEM \\
by \\
James Ulrich \\
\end{center}

\smallskip

\noindent Advisor: Professor Michael Anshel. \\

	Motivated by a desire to test the security of the pubic key exchange
protocol of I. Anshel, M. Anshel, and D. Goldfeld, (``An Algebraic
Method for Public-Key Cryptography'', Mathematical Research Letters,
vol. 6, pp. 1-5, 1999), we study algorithmic approaches to the
simple commutator decision and promise problems (SCDP/SCPP) for the
braid groups $B_n$. We take as our point of departure a seminal
paper of O. Ore, (``Some Remarks on Commutators'', Proceedings of
the American Mathematical Society, Vol. 2, No. 2, pp.307-314, 1951),
which studies the SCPP for the symmetric groups.

Our results build on the work of H. Cejtin and I. Rivin, (``A
Property of Alternating Groups'', arXiv:math.GR/0303036). We
extract, from their proof that any element of the alternating
subgroup of $S_n$ can be written as a product of two $n$-cycles, an
explicit algorithm for solving the SCPP for $S_n$. We define a model
of computation with respect to which the algorithm executes in time
$O(n^2)$.

We then extend the algorithm to a subset of permutation braids of
the braid groups $B_n$, to show that any element of the commutator
subgroup $[B_n,B_n]$ may be efficiently written as the product of a
pure braid and a simple commutator of permutation braids. We use
this result to define a probabilistic approach to the SCDP/SCPP,
posing for future research the question of whether such an algorithm
may be made efficient with respect to a measure of complexity such
as that defined in a work of I. Kapovich, A.  Myasnikov, P. Schupp,
V. Shpilrain (``Average-Case Complexity and Decision Problems in
Group Theory'', Advances in Math. vol. 190, pp. 343-359, 2005).

\end{abstract}

\pagebreak

\begin{center}
\textsc{Acknowledgements}
\end{center}
I wish to thank my advisor, Professor Michael Anshel, for his steady guidance over the course of my graduate career. 
It is only slightly less a tautology than the statement $1=1$ to say that without him, this work would not have been
possible. I also wish to thank the other members of my defense committee, Professors Cormac
O'Sullivan, Burton Randol, and Vladimir Shpilrain, for their time and helpful advice.  Thanks
are also due Professors Joan Birman, Edgar Feldman, Minhyong Kim, Roman Kossak,  Henry Pinkham, Dennis Sullivan, Lucien Szpiro, and Alphonse Vasquez for their generous assistance at various key points of my studies. Thanks too are due my colleagues  Tara Brendle, Arjune Budhram, Hessam Hamidi-Tehrani, and Brendan Owens for
their support, educational and otherwise,  for lo these many years.  Of course, I must also thank Jocelyn, my love and partner,
for her general willingness to put up with me, as well as my mother and father, Mary
Louise and David, for scraping together between them just enough math DNA to give me a fighting chance.

\vspace*{4.0in}

\pagebreak

\tableofcontents

\pagebreak

\listoffigures

\pagebreak

\pagestyle{myheadings}
\pagenumbering{arabic}

\section{Introduction: Ore's commutator problem.}

In the August 2004 issue of the Notices of the American Mathematical
Society \cite{Aschbacher}, Michael Aschbacher reported on the state of
the Gorenstein, Lyons, and Solomon program, begun in the 1980s, to
establish a formal, cohesively written proof of the classification theorem
for simple finite groups. The theorem states that all
finite simple groups fall into one of the following classes:
groups of prime order, alternating groups, groups of Lie type
(that is, having a representation involving automorphisms of a vector
space over a finite field), or one of 26 ``sporadic'' groups (that
is, exceptions to the preceding). The classification theorem is
central to the study of finite groups $G$, since the simple
factors of a composition series for G, in Aschenbacher's words,
``exert a lot of control over the gross structure of $G$.'' (Recall
that a composition series for $G$ is a sequence of normal
subgroups $1 = G_0 \lhd G_1 \lhd \cdots \lhd G_n = G$, where each
$G_i$ is simple -- that is, contains no normal proper subgroups).

Accordingly, a conjecture given by Oystein Ore in his seminal 1951
paper ``Some remarks on commutators''  \cite{Ore} has been of interest
to researchers concerned with the classification problem. In that
paper, Ore studies the \emph{symmetric group} $S_n$ (which we recall
is the group of permutations of a set of $n$ elements), and its
alternating and derived subgroups.  (Recall that the
\emph{alternating subgroup} $A_n \subset S_n$ is the subgroup of
permutations that can be written as products of an even number of
transpositions -- that is,  swaps -- of adjacent elements. Recall
also that the \emph{derived group} or \emph{commutator subgroup}
$S_n^\prime \subset S_n$ is the group  generated by the
\emph{simple commutators} of $S_n$,  which are elements  of the form
 $[x,y] := xyx^{-1}y^{-1}$, for $x, y$ in $S_n$).  In general, elements of the commutator
 subgroup of a given group are not themselves simple commutators (see \cite{Cassidy}, \cite{Isaacs}).  
 Ore conjectures in the paper that \emph{every} element of a simple group $G$ of finite order is in 
 fact a simple commutator of elements of $G$. A key result of his paper is:
 \begin{prop}[Ore, \cite{Ore}, theorem 7)] \label{prop:ore}
\emph{For $n \geq 5$, every element of the alternating group $A_n$ is a simple commutator
of elements of $A_n$}.
\end{prop}

The authors Ellers and Gordeev, in ``On the Conjectures of J. Thompson and O. Ore"  \cite{Ellers},
note that  a stronger conjecture is attributed to R. C. Thompson: \emph{every finite simple group G
contains a conjugacy class $C$ such that $C^2 = G$}, which implies Ore's conjecture.Ê (An explanation of why
Thompson's conjecture  implies Ore's  is given in \cite{Arad}).  The authors describe many examples of groups
for which the Thompson conjecture is known to be true, including the projective special linear group
$ \text{PSL}_n(K)$ for $K$ a finite field, and show that the conjecture holds for all groups of
Lie type over finite fields containing
more than $8$ elements. \\
\\
The work of Cejtin and Rivin \cite{Rivin} is of particular interest to us, as it asserts the following:

\begin{prop}[Cejtin and Rivin, \cite{Rivin}]  \label{prop:rivin}
\emph{There is an efficient algorithm to write every element $\sigma$ of $A_n$ as a product of two $n$-cycles}.
\end{prop}
From this, the authors show that there is an efficient algorithm to solve the \emph{simple commutator promise problem} (SCPP) for
$S_n$.  In general, given an arbitrary group $G$,  and $g \in G$ guaranteed to be a simple commutator, the SCPP for $G$ asks
for explicit $x,y \in G$ such that $g=[x,y]$.

We will be concerned here with the symmetric groups $S_n$ and also with the braid groups $B_n$,
defined below. The braid groups play a central role in knot theory and the topology of $3$ and $4$-dimensional manifolds 
\cite{Adams}  \cite{Birman} \cite{Kirby} \cite{Witten} \cite{Gompf} \cite{Ohtsuki}. They also play
a significant role in the public key exhange protocol of Anshel, Anshel,
and Goldfeld \cite{Anshel}. Hence finding an efficient method of solving the SCPP for $B_n$ is
an area of active research, as is finding efficient methods for solving
the related \emph{conjugacy search problem}: given elements $u,w \in B_n$, find $x \in B_n$
such that $u = xwx^{-1}$.

In what follows below, we will examine the SCPP for the symmetric groups $S_n$ and the braid groups $B_n$. Where
the braid groups are concerned, we will restrict the problem to those
elements $w$ of a braid group $B_n$ that are simple commutators $w = [x,y]$ of permutation braids  $x$, $y$,
where a permutation braid is a positive braid such that no two strands cross twice.  Any element of a braid
group $B_n$ is a product of permutation braids and their inverses, and for a given $B_n$, there
is a bijective map between the set of permutation braids of $B_n$ and $S_n$.  We will examine ways
in which the Cejtin-Rivin algorithm can be used to address the SCPP for $B_n$, restricted to simple commutators
of permutation braids.

To make the discussion of computational complexity somewhat more rigorous, we will first define our model of a classical Turing machine, along with our versions of the notions of alphabets, words, and languages, in order to map between algorithms expressed in terms of Turing machines and those expressed through algebraic and symbolic operations (i.e. ``psuedo-code''). This will allow us to
define our notions of computational complexity. We will then describe the Cejtin-Rivin algorithm for the simple commutator
promise problem for $S_n$, and discuss its complexity. We will present an explicit program, in terms of algebraic and symbolic
operations, to implement the algorithm. Following this, we will provide definitions of the relevant concepts of braid groups.
We will explore the extension of the Cejtin-Rivin algorithm to simple commutators of elements of permutation braids. Finally,
we will describe possible avenues for future research.

\pagebreak

\section{Classical Turing machines and computational complexity}

In this section we describe the classical Turing machine, a model of computing that will allow us to discuss
algorithmic complexity.

\subsection{Classical Turing machines}

In order to speak about the complexity of computational problems,
we need to have some sort of computer in mind. So we define the
notion of a Turing machine, which is a simple but powerful model
of a computer, generally attributed to Alan Turing (1912-1954). We
use the formulation given by \cite{Lewis}.

\begin{defn}\label{def:alphabet}  For our purposes, an \emph{alphabet} $\Sigma$ will denote
 a subset of the set whose elements are the  upper and lower-case letters of the English alphabet, 
 the digits $0$ through $9$, the symbols $\sqcup$ and $\triangleright$,
 and the standard English punctuation symbols.
 A \emph{language} $L$ is then a subset of $\Sigma^*$, the set of all finite strings of
symbols from $\Sigma$. So, for example, if $\Sigma = \{ 0,1 \}$
then we might have $L \subseteq \Sigma^*, L =
\{0,10,100,110,1000,1010,\cdots\}$, the set of all even numbers,
expressed in binary form.
\end{defn}

\begin{defn}\label{def:turing_machine}
A \emph{classical, single-tape, deterministic Turing machine (CDTM)} consists of: \\
\\
(i) a quintuple $ (K,\Sigma,\delta,s, H ) $ where $K$ is
a finite set of \emph{states}; $\Sigma$  is an alphabet
containing the \emph{blank symbol} $\sqcup$ and the \emph{left
end symbol} $\triangleright$; $s_i \in K$ is the \emph{initial state}; $s_f \in K$
is the \emph{halting state}; $\delta$
is the \emph{transition function} from $ K \times \Sigma $
to $ K \times ( \Sigma \cup \{ \leftarrow, \rightarrow \} ) $ such
that for all $q \in K $, if $\delta(q,\triangleright) = (p,b)$
then $b = \rightarrow$ and such that for all $q \in K $ and $a
\in \Sigma$, if
$\delta(q,a) = (p,b)$ then $b \neq \triangleright$. \\

\noindent (ii) a unidirectionally infinite (to the right)
\emph{tape} consisting of squares, each containing one symbol
from $\Sigma$, a finite number of which are not the blank symbol.
The \emph{input} $w \in \Sigma^*$ consists of symbols from
$\Sigma - \{ \triangleright, \sqcup \}$, and follows the
$\triangleright$ at the left end of the tape; the
first blank symbol denotes the end of input. \\

\noindent (iii) a \emph{tape head} that can read from and write to
the tape, and move left and right along the tape. 
\end{defn}

It is held that what can be computed on any given existing (classical) computer may
also be computed on this version of a Turing machine; this assertion is
known as the \emph{Church-Turing thesis} (see \cite{Nielsen} p. 125).
The machine  works as follows. The machine starts in the initial state $s_i$,
with the tape head positioned over the square immediately to the right of
the left end symbol. The tape reads the symbol $b$ at that square.
The transition function $\delta(s_i, b)$ then yields a tuplet
$(k,b^\prime)$. If $b^\prime = \rightarrow$, the tape head moves
right one square; if $b^\prime = \leftarrow$, the tape head moves
left one square (if not already at the left end of the tape);
otherwise, the tape head replaces the symbol in the current square
with $b^\prime$. The machine then enters the state $k$. If $k = s_f$,
the machine halts, and the contents of the tape at that time is said to be
the $\emph{output}$ of the machine. 

\begin{ex}[addition of two $n$ digit binary numbers] \label{ex:turing_addition}
Let $\Sigma$ be the alphabet $\{0,1\}$ and let $L$ be the set
of all pairs of symbols from $\Sigma$; that is,
$L = \{(0,0),(0,1),(1,0),(1,1) \}$.  Then here is a description
of the Turing Machine that accepts a pair from $L$, and
adds the two elements of the pair together, outputting the result.
We assume that the input $w$ consists of the sequence $a b$
where $a$ is the first element of the pair, and $b$ the second element.
For convenience, we assume that the machine output will begin on
the third square of the tape. \\

\indent \emph{Turing Machine to add two $1$-digit binary numbers:} \\
\indent  \indent \emph{input}:  Two binary digits in successive tape squares.\\
\indent  \indent \emph{states}:  $K = \{s_i,s_f, s_{c}, s_{nc}, s_{0b}, s_{1b}  \}$ \\
\indent \indent \indent \indent \emph{transition function $\delta$}:\\
\begin{equation} \nonumber
\begin{array}{lll}
d(s_i,1) =   (s_c, \rightarrow) & d(s_c,1) =  (s_1, \rightarrow) & d(s_1,\sqcup) =  (s_1,1)  \\
d(s_1,1) =  (s_0,\rightarrow)  &d(s_0,\sqcup) =  (s_f,0)  & d(s_i,1) =  (s_{c},\rightarrow)  \\
d(s_{c},0) =  (s_{1b},\rightarrow) & d(s_{1b},\sqcup) =  (s_f,1) & d(s_i,0)  =  (s_{nc},\rightarrow)  \\
d(s_{nc},1)  =  (s_{1b},\rightarrow)  & d(s_i,0)  =  (s_{nc},\rightarrow)  &  d(s_{nc},0)  =  (s_{0b},\rightarrow) \\
d(s_{0b},\sqcup)  =  (s_f,0) & &  \\
\end{array}
\end{equation}
\\
\end{ex}
%\vspace{-0.5in}
The machine examines the contents of the first square; if the content is a $1$, it enters the state
$s_c$; otherwise it enters the state $s_{nc}$. It then advances to the next square. If the state is
$s_c$ and the content of the second square  is $1$, it advances to the next square, outputs $10$, and
terminates. Otherwise if the state is $s_c$ and the contents of the second square is $0$, it advances
one square, outputs $1$, and terminates. It behaves in similar fashion for the case in which the content
of the first square is a $0$.

More generally, one can describe a single-tape Turing machine that adds two $n$-digit (binary) numbers, in part
as follows:\\
\\
\indent \emph{Turing Machine to add two $n$-digit binary numbers:} \\
\indent \indent \indent \indent \emph{input}:  Two $n$-digit binary numbers $a$ and $b$, each terminated by
a blank, and printed on the input tape left to write (least significant digit of each number on the left).\\
\indent \indent \indent \indent \emph{states}:  $K =   s_i, a_1, b_1, c_1, c_{1c}, r, r_c, n_c, \cdots $ \\
\indent \indent \indent \indent \emph{transition function $\delta$}:\\

\noindent \emph{state/transitions for the for the case in which we have no carry,
and both the current digit of $a$ and the current digit of $b$ are $1$}:
\begin{equation} \nonumber
\begin{array}{lll}
d(n,1) =  (X, a_1)  & d(s_i,1) =  (X, a_1) & d(a_1,X) =  (\rightarrow, a_1) \\
d(a_1,! \sqcup) =  (\rightarrow, a_1) & d(a_1,\sqcup) = (\rightarrow, b_1) & d(b_1,Y) =  (\rightarrow, b_1) \\
d(b_1,1) =  (Y, c_{1c}) & d(c_{1c}, ! \sqcup) =  (\rightarrow, c_{1c}) & d(c_{1c}, \sqcup) =  (1,r_{c}) \\
d(r_c, !X) =  (\leftarrow 1,r_{c}) & d(r_c, X) =  (\rightarrow1,n_c) \\
\end{array}
\end{equation}

Here $!S$ denotes ``any symbol other than S.'' The machine starts with the tape head at the leftmost digit of $a$, and reads the digit. If the digit is $1$, it enters the state $a_1$ (to record that it read a $1$ from $a$), marks the square with an $X$, and 
moves the tape head to the right.  If it encounters any symbol other than a blank, it remains in state $a_1$, and continues 
to move to the right. Otherwise, it enters the state $b_1$ (to record that it is now processing digits of $b$, having read a $1$ 
from $a$), and moves the tape head to the right. It continues to move to the right until it stops at the  next digit of $b$ to be processed (this will be the  first non-$Y$ digit). If that digit is $1$, it marks the square with a $Y$, enters the state $c_{1c}$, (to record that it should output a $1$, and that a carry occurred)(, and moves the tapehead to the next output square. There it records a $1$, and enters the state $r_c$ (to record that a carry has occurred, and that it must now reverse to the left to find the next digit of $a$ to process). The sets of states corresponding to the other scenarios (current digit of $a$ is $0$, current digit of $b$ is $1$, carry or no carry), as
well as end of input logic, follow similarly. 
 
\subsection{Concering Computational Complexity}

An inspection of the example of the Turing machine given for adding two $1$ digit numbers reveals that it will require at most $5$
invocations of the transition function $\delta$. The more general machine that adds two $n$-digit binary numbers needs
to invoke the transition function on the order of $2n$ times, for each pair of digits it processes (one from each number
to be added). This is because the single tape head must move back and forth across the numbers as it adds them. It follows that 
the single tape, classical, determinstic Turing machine $M$ described above will add two $n$-digit binary integers in a number of invocations $Cn^2 + B$ of the transition function $\delta$ of $M$, for constants $C$ and $B$. \\
\\
In general, given a Turing machine $M$ designed to compute some problem, one asks for an upper bound $f(l)$
on the number of required invocations of the $\delta$ function of  $M$,  terms of the length $l$ of the input. Such an upper
bound provides a rough measure of the \emph{time complexity} of the problem. Hence, we would say that there is
a \emph{quadratic-time} algorithm to compute the addition of two binary integers, or that addition of integers
is ``quadratic time,'' with respect to the Turing machine described above. (We equate one invocation of the transition function with one ``clock tick.'')

\begin{defn}[order of complexity]\label{def:complexity}
One says that an algorithm executes in time $O(f(l))$, where $f$ is a given function of the length $l$ of the input, 
if there exists \emph{some} Turing machine $M$, such that there exist constants $c$ and
$l_0$ such that, for all values of $l$ greater than $l_0$, the number $N$ of invocations of the transitition
function of $M$ required by the algorithm satisifes $N \leq cf(l)$. If there exists a polynomial function $f(l)$
for which this is true, we say the algorithm executes in ``polynomial time'', and
we say that it is a member of the (time) complexity class $\mathbf{P}$ of polynomial-time algorithms. Such algorithms
are considered ``tractable'' -- that is, amenable to computation. By contrast, problems having greater than
polynomial-time complexity are considered ``intractable''. 

\end{defn} 
\noindent 

\begin{defn}[complexity in terms of higher-level programming languages]  \label{def:equiv_comp}
In the following, we will consider an algorithm for the SCPP for $S_n$, given by Cejtin and Rivin in \cite{Rivin}, to which the authors have assigned the time complexity $O(n)$.  As explained in \cite{Nielsen}, pp.135-142, we will adopt the hypothesis that if we can
specify an algorithm in a high-level programming language such as C, C++, or Java, or in a pseudo-code equivalent (defined below),
then there is a  multi-tape, classical, deterministic Turing machine, following that given in \cite{Domanksi}, such that it can act as an ``interpreter'' for the algorithm.  That is, the machine $M$ can accept the pseudo-code algorithm $A$ and a given input $I$ to the algorithm, and apply the algorithm $A$ to the input $I$. We further suppose that it can do so  in a way such that, for each execution of a pseudo-code statement,  the machine $M$ invokes its transition function once and only once. Then we say that the algorithm is $O(f(n))$ if, for any algorithm input $I$ of length $n$, the machine will complete processing  in $O(f(n))$ transition invocations. We will say that $A$ is efficient if the corresponding number of transition invocations $f(n)$ is a polynomial in $n$. 
\end{defn}
  
\pagebreak

\section{The simple commutator decision problem for $S_n$.}

In this section we describe the Cejtin-Rivin algorithm for the simple commutator promise problem 
for the symmetric groups $S_n$, provide an implementation of the algorithm in pseudo-code,
and discuss its complexity.

\subsection{Preliminaries concerning $S_n$}

\begin{rem}[a definition of $S_n$]\label{rem:def_of_$S_n$}

Recall again some basic facts concerning the symmetric group $S_n$, which is is the group of permutations
of the set  $\{1,2,\cdots,n\}$.  That is, a permutation $\sigma \in S_n$ is a $1-1$ map from the set
$\{1,2,\cdots,n\}$ to itself. One may imagine the elements being the points
$\{(0,1),(0,2),\cdots,(0,n)\} \in \mathbb{R}^2$, with point $(0,i)$ initially labelled by the label
$l_{(0,i)}= i$.  A permutation $\sigma$ is then specified by a re-assignment  of the labels associated to the points:
$l_{(0,i)} \rightarrow \sigma((0,i))$.  For example, consider the element $\sigma \in S_3$ given by
$(l_{(0,1)} = 1, l_{(0,2)} = 2, l_{(0,3)}= 3) \rightarrow (l_{(0,1)} = 2, l_{(0,2)} = 1, l_{(0,3)} = 3)$. This is the element
that interchanges the labels associated to the points $(0,1)$ and $(0,2)$.

The permutations form a group under composition, the group being generated by the standard generators $\tau_i$,
where $\tau_i$ denotes the interchange of the labels associated to the elements in the $(0,i)$ and $(0,i+1)$-th positions.
(From now on, we will denote the element $(0,i)$ simply by the integer $i$, and speak of ``permuting elements $i$ and $j$,''
though we really mean to permute the labels associated to $(0,i)$ and $(0,j)$. So in our above example, we would
write that $\sigma(1) = 2, \sigma(2) = 1$, and $\sigma(3) = 3$.)   When mutiplying permutations, we read
from right to left, hence $(1,2)\circ (2,3)$ means ``first exchange the elements  $2$ and $3$, then
exchange the elements $1$ and $2$.''

\end{rem}

\begin{rem}[alternating subgroup $A_n$ of $S_n$]\label{rem:alternating_subgroup}

The alternating subgroup $A_n$ of
$S_n$ is the group of permutations that may be expressed as a product of an even number of adjacent transpositions.
Such permutations are called \emph{even}.
\end{rem}

\begin{rem}[cyclic decomposition of elements of $S_n$]\label{rem:cyclic_decomposition}
One may record a permutation $\sigma \in S_n$ in terms of its \emph{cyclic decomposition},
a cycle being a circular sequence $i \rightarrow j \rightarrow k \rightarrow \cdots \rightarrow i$
of elements exchanged by the permutation.  One typically elides the $1$-cycles. So, for example, the
permutation $\sigma \in S_5$ taking the elements $1,2,3,4, 5$ to $3,4,1,2,5$ could be
written $\sigma = (13)(24)(5) = (13)(24)$. The order in which the cycles appear does not
change the  corresponding permutation, nor do cylic permutations of the sequence in which the integers of
any given cycle appear. (So $(13)(24)$ and $(42)(31)$ denote the same permutation).
\end{rem}

\begin{rem}[order of a cycle]\label{rem:cycle_order}

 The number of elements appearing in a cycle is known as the \emph{order} of the cycle; a specification of
the number of cycles of each order appearing in $\sigma$, counted by multiplicity, gives the
\emph{type} of the cyclic decomposition of of $\sigma$. So, in our example, $\sigma$ consists
of two $2$-cycles (and implicitly, a single $1$-cyle). A permutation $\pi$ will be conjugate to a given permutation $\sigma$ if and only if
$\pi$ and $\sigma$ are of same  cyclic decomposition type.  An element of $A_n$ will consist of cycles of
odd order, and an even number of cycles of even order.
\end{rem}

\begin{rem}[standard representation of $S_n$]\label{rem:standard_rep}
 Finally, we obtain a faithful linear representation $\rho$ of $\Sigma$ as the matrix group generated
 by the images of the standard generators $\tau_i$ under
$\rho$, where $\rho(\tau_i)$ is the $n \times n$ matrix given by interchanging the $i$ and $i+1$-th columns
of the indentity matrix of $GL(n,\Real)$.
\end{rem}

The information in these remarks may be found in  \cite{Sagan} and \cite{Ore}.

\subsection{The Cejtin-Rivin algorithm for the SCPP for $S_n$}

The paper ``A Property of Alternating Groups,'' co-authored by Henry
Cejtin and Igor Rivin  (\cite{Rivin}) shows [theorem 1, page 1] that
there is an efficient algorithm to write any even permutation $\sigma$
in $S_n$ as a product of two $n$-cycles $p_1$ and $p_2$; that is,
$\sigma = p_1 p_2$. Since $p_1$and $p_2$ are both $n$-cycles, it
follows from remark \ref{rem:cycle_order} that $p_2 = \nu p_1
\nu^{-1}$ for some $\nu \in S_n$.  Moreover, if $p_1$ is an
$n$-cycle, then so is $p_1^{-1}$: if $p_1 = (i_1 i_2 \cdots i_{n-1}
i_n)$, then $p_1^{-1} = (i_n i_{n-1} \cdots i_2 i_1)$. Hence, $p_2$
and $p_1^{-1}$ are both $n$-cycles, too. Thus $p_2 = \tau p_1^{-1} \tau^{-1}$
for some $\tau \in S_n$. So we have $\sigma = [p_1,\tau]$, for some $\tau \in S_n$. 

Moreover, we can readily identify such a $\tau$, since if 
$p_1^{-1} = (i_1 i_2 \cdots i_{n-1})$ and  $p_2 = (j_1 j_2 \cdots j_n)$, 
then $\,\tau$ is the unique permutation such that 
$\tau(i_1) = j_1, \tau(i_2) = j_2, \cdots ,\tau(i_n) = j_n$ (see \cite{Sagan}, p. 3). \\
\\
 We provide a brief summary of the (constructive) proof of the theorem as
given in the paper.  We will then provide a detailed implementation of a functional equivalent of the
algorithm constructed in the proof; we will express the algorithm in pseudo-code, as defined below.
A complexity analysis of the implementation will then yield a concrete verification that the algorithm 
posited by the theorem is indeed efficient, providing an efficient solution for the $SCPP$ for $S_n$. 

Note: the authors place the algorithm in the complexity class $O(n)$, but do not specify a model of computation. 
Moreover, the algorithm is not given explicitly as psuedo-code, but implicitly, in the course of proving constructively,
by induction, the \emph{existence}  of a method for writing any element $\sigma \in A_n \subset S_n$ as a product of
two $n$-cycles.

We show below that, using the definition (\ref{def:equiv_comp}), our pseudo-code implementation will 
have complexity $O(n^2)$, for  input 
a word $w$ in the standard generators of $S_n$, having fixed length $k$, for variable $n$, and that it
will execute in time $O(k)$, for fixed $n$ and variable $k$ .\\

\subsubsection{Review of Cejtin-Rivin proof}
Now, in the paper, the authors first show (lemma $3$, p.2), that given a permutation $p$ in $S_n$, such that:\\
\\
1) it is the product of two permutations $s_1$ and $s_2$ \\
\\
2) $s_1$ acts non-trivially only on ${1,\cdots, t}$, and $s_2$ acts non-trivially only on ${t+1,\cdots,t+s}$ where $n = t + s$\\
\\
3) $s_1 = p_{11} p_{12}$ where $p_{11}$ and $p_{12}$ are $t$-cycles\\
\\
4) $s_2 = p_{21} p_{22}$ where $p_{21}$ and $p_{22}$ are $s$-cycles\\
\\
then $p = s_1 s_2 = (p_{11} p_{21} v)(v p_{12} p_{22})$ is a product of two $(s+t)$-cycles, where $v = (t, t+s)$.\\
\\
The authors next prove: \\
\\
Lemma $4$: if $n$ is odd, then any $n$-cycle $\sigma$ is a product of two $n$-cycles $\rho_1$ and $\rho_2$. \\
\\
Lemma $5$: if $n = 4m$, and $\sigma \in S_n$ is a product of two disjoint $2m$-cycles, then $\sigma$ is a product
of two $n$-cycles $\rho_1$ and $\rho_2$.\\
\\
Lemma $6$: if $n = 2s + 2t$, where $s < t$, and if $\sigma \in S_n$ is a product of a $2s$-cycle and a $2t$-cycle, such
that the cycles are disjoint, then $\sigma$ is a product of two $n$-cycles $\rho_1$ and $\rho_2$.\\
\\
The proofs of lemmas $4,5,6$ each provide explicit cyclic decompositions of the $n$-cycles $\rho_1$ and $\rho_2$.\\
\\
In the proof of the theorem itself, it is first noted that the theorem holds for $n=1$. The proof proceeds by induction.
Given an element $\sigma \in A_{n+1}$, the alternating group of $S_{n+1}$, if the element consists of exactly one (odd) cycle,
or if it is a product of exactly two (necessarily even) cycles, then the proof follows from the
constructions given in  lemmas $4$,$5$,  and $6$. If, on the other hand,
$\sigma$ is a product of  more than two cycles, it follows  that it can be written  as a product of two disjoint
even permutations $c_1$ and $c_2$, or else as the product where $c_1$ acts non-trivially only on ${1,\cdots,t}$ and $c_2$ acts
non-trivially only on ${t+1,\cdots t+s}$  where $t+s = n+1$.  The inductive hypothesis and lemma $3$
then combine to say that $\sigma$ is a product of two ${n+1}$ cycles: the inductive hypothesis says that $c_1$ is a product
of two $t$-cycles $p_{11}$ and $p_{12}$ and $c_2$ is a product of two $s$-cyles $p_{21}$ and $p_{22}$; lemma $3$ says that
$\sigma$ is then a product of two $s+t$-cycles.

One must still obtain explicit expressions for the $p_{ij}$ in this latter case. The authors show that  $c_1$
can always be written as a product of exaclty one or two cycles, from which the expressions for $p_{11}$ and $p_{12}$ follow
directly from lemmas $4,5$, and $6$. To obtain the expressions for $p_{21}$ and $p_{22}$, one must recursively apply the argument
of the proof to $\sigma^\prime = c_2$; clearly no more than $n-3$ recursions will be required, since any $3$-cycle can
be written as a product of no more than two cycles.   What remains to be shown is that the entire process can be coded into
an algorithm that can execute on the Turing machine described above in an efficient manner.\\

As promnised above,  we next present an explict implementation of the Cejtin-Rivin algorithm, using ``pseudo-code'' that 
can be translated both to the language of Turing machines as described above, and to modern computer programming 
languages such as C++. We will then assess the complexity of the algorithm to see how far it varies from the $O(n)$ bounds 
set by Cejtin and Rivin. It is clear from the discussion above that any reasonable implementation of the Cejtin-Rivin algorithm 
will fall in the complexity class  $P$.

\subsection{Implementation of Cejtin-Rivin in pseudo-code}

We here present an explicit implementation of the Cejtin-Rivin algorithm described above. The algorithm is expressed in
pseudo-code. It is difficult to define pseudo-code in mathematically precise terms.  Here, we will define psuedo-code
to be a set of human-readable instructions, such that they would provide, to a human trained in the art of computer programming,
all the information necessary to create a computer program, written in the computer language C++, to perform an algorithm such as that developed by Cejtin-Rivin, as described above.  (More concretely, we can say that pseudo-code must have the property
that every statement maps in an ``obvious'' way to a programming language such as C, C++, or Java). We assume in our psuedo-code that we are free to utilize the following types of programming operations, each of which has a direct 
implementation in terms of the Turing machine defined above, and also in a higher-level programming language such as C.

\subsubsection{Types of operations comprising psuedo-code to be used below}\label{sub:operations}

The types of operations to be used in the psuedo-code below are:\\
\\
\emph{Operation Type $1$}:  read a string symbols from an input device. \\
\\
\emph{Operation Type $2$}: write a string of symbols to an output device. \\
\\
\emph{Operation Type $3$}: declare a program variable. Program variables are integers, or lists thereof.
Lists of integers are said to be \emph{arrays} and have the property that the integer in the $i$-th position
of the array $s$ (arranging the list horizontally, and counting left to right,  starting with $0$), may be written or  
read via the notation $s(i)$. (Observe that we can represent a matrix as an array of arrays).\\
\\
\emph{Operation Type $4$}:  assign a value to a variable. We may assign integral values to integers (using the notation
$i = n$ to assign the value $n$ to integer $i$) , and lists of integral
values to arrays of integers.  In a Turing machine, a program variable may be represented by the reservation of a fixed amount
of space beginning at a certain location on the storage tape. We assume that all program variables occupy appropriate storage
to make assignment statements valid. \\
\\
\emph{Operation Type $5$}: perform a logical operation or algebraic operation, such as testing whether a program
variable $i$ is zero, or comparing two integers to see if they are equal, or adding one integer to another. The equality
 operation is denoted $i == j$; the ``test for zero'' operation will be denoted  $!i$.  \\
\\
\emph{Operation Type $6$}: iterate a series of instructions for some count. That is, we have access to iterative control structures
such as ``while condition $A$ is true, do sequence of statements $B$'' or
``if condition $A$ is true, then do sequence of statements $B$'' or ``for integer $x$ = 0, while  $x <  N$, do sequence
of statements $B$, then increment $x$ by $1$''. Where we have a sequence of instructions  that are to be implemented as
part of a control structure, we indent the instructions underneath the instruction defining the control structure condidtion.
Numbering each statement within an algorithm, we also assume the existence of the infamous ``Go To (line number)'' control
operation. \\
\\
\emph{Operation Type $7$}: call a subalgorithm, or ``function''. This allows us to organize the transistions
of our Turing machine into functions, which are self-contained subalgorithms accepting input and producing output.
Functions may invoke other functions. Lists of variables of varying types (called ``arguments'') may be passed from
the invoking  function to the function invoked, and lists of variables of varying types may also
be passed from the invoked function back to the invoking function.

We use ``call F(arg1,arg2,...,argN)'' to denote the invocation of function $F$, passing input ``arg1,arg2,..,argN''.
We use, within a function, ``return(arg1,arg2,...,argN)'' to denote return to the invoking function, passing back
``arg1,arg2,...,argN''.  A function may invoke itself. 

All these operations are directly expressible in a wide variety of computer programming languages.  We also allow comments
to be interspersed with the psuedo-code;  we place comment text between ``/*'' and  `` */'' on the line preceding the
psuedo-code line to which the comment applies.  The comments are not regarded as part of the algorithm.

\subsection{An aside: complexity of the preceding constructs}

\begin{lem}\label{lem:psuedo_complexity}
The constructs $1-7$ of the preceding section all execute on a classical, single-tape, deterministic 
Turing machine in an number of operations that is of order $O(n^2)$, where
$n$ is the size (in tape cells) of the input to the operation. \\
\end{lem}

\noindent \emph{proof}:
All complexity bounds are given with respect to the single-tape, deterministic Turing machine described above. For all operations, we assume a fixed number of tape squares between the tape head at the beginning of the operation, and the nearest tape location at which data is to be read or written (i.e. we do not count any  intervening``extraneous'' tape data that we must navigate
around due to our machine having just one tape).

Operations of type $1$ and $2$ are $O(n)$ by construction, and  definition \ref{def:turing_machine}.

Operations of type $3$ amount to reserving a certain blank portion of the storage tape of \ref{def:turing_machine}, and
writing an initial value to that section of the tape. This operation is $O(n)$ because
operation $2$ is.  Operations of type  $4$  involve overwriting a section of the storage
tape reserved by an operation of type $3$.   Operations of type $5$ amount to comparison of an integer to determine
if all digits of the integer are $0$ (``test for zero''),  or addition or subtraction of two integers, or to
subtraction of two integers, after which the sign of the result is checked (comparison, equality testing).
Note that subtraction may be performed in the same manner as addition, by first converting the subtrahend to its negative. By example \ref{ex:turing_addition}, addition is an $O(n^2)$ operation; negatation of a binary number simply entails the inversion
of all $0$s and $1$s, followed by the addition of the number $1$ (``two's complement'' representation of a number).

Operations of type $6$ involve the repeated invocation of a given sequence of transition functions, where the number
 of invocations is controlled by a counter which is decremented or incremented after each invocation. The invocations are terminated when the counter reaches a set value. Since initialization of an integer, and the subtraction or addition of two integers, are at most $O(n^2)$ operations, the act of iteration is $O(n^2)$ (without consideration of the thing iterated). Here $n$ is
the storage size of the counter in tape cells.

Operations of type $7$ amount to copying a series of symbols 
to a blank section of the storage tape, invoking a set sequence of
transition functions which produce output on another section of the
storage tape, and then reading the result of that output. Hence
function invocation is simply an operation of type $2$, followed by
a series of operations of type $1$ (without regard to the complexity of the
invoked function). Hence an operation of type $7$ is of order $O(n)$.  \indent  $\Box$

As explained in definition (\ref{def:equiv_comp}), in what follows we will regard each operation of types $1-7$ as constituting
one clock tick.

\pagebreak
\subsection{Algorithm:  Expression of an element of $A_n$ as a simple commutator of elements of $S_n$}

We here present, using the constructs of the psuedo-code language
described above, an implementation of the Cejtin-Rivin algorithm for
the SCPP for $S_n$ \cite{Ulrich}.\\
\\

\noindent \emph{PROBLEM:  Express $n \in A_n$ as a simple commutator.} \\
\\
\noindent INPUT: A positive integer $n$, and an element $\sigma$ of the 
alternating subgroup $A_n$ of $S_n$, expressed as a word $w(\sigma)$ 
in the standard generators $\tau_i$ of $S_n$, where $\tau_j$ is denoted by the integer $j$.\\
\\
\noindent OUTPUT: The element $\sigma$ expressed as a simple
commutator of elements $\rho_1$ and $\rho_2$ of $S_n$; that is, an
expression of $\sigma$ as  $\rho_1 \rho_2 \rho_1^{-1} \rho_2^{-1}$.
The elements $\rho_1$ and $\rho_2$ will each be expressed as words in
the standard generators of $S_n$; with the integer $i$ denoting
the generator $\tau_i$. Each word will be separated by a $0$, and the
the entire output will be terminated by $00$.  The output should be
read left to right (and applied right to left).\\
\\
ALGORITHM:\\
\\
\indent 1. We have that the input $w$ is of the form $w = i_1 i_2
i_3 \cdots i_k$ where each $i_j$ is an integer between $1$ and
$n-1$. 

\indent 2. Perform subalgorithm $A$, passing $(w,n)$  as input. Let the
output of $A$ be denoted by $M$. Then $M$ is an $n \times n$  matrix
that is the standard representation of the permutation given by $w$. 

\indent 3. Perform subalgorithm $B$, passing $(M,n)$ as input. Let the
output of $B$ be denoted by $d$. Then $d$ is a sequence of integers
representing the cyclic decomposition of the permutation represented
by $M$. 

\indent 4. Perform subalgorithm $C$, passing $(d,n)$ as input. Let
the output of $C$ be denoted by the pair $(c_1,c_2)$, where $c_1$
and $c_2$ are both products of cycles, with each cycle being given
by a sequence of integers. (For the exact format of $c_1$ and $c_2$,
see the description of the output of subalgorithm $B$).

\indent 5. We have that $c_1$ and $c_2$ are each a product of
cycles, expressed in the form $i_1 i_2 \cdots i_{j-1} i_j 0 i_{j+1}
\cdots i_k 0 \cdots 0 i_l \cdots i_m 00$. Let $S$ denote the matrix
output by subalgorithm $E$, invoked with input $c_1$. Perform the
same steps for $c_2$ to obtain a matrix $T$. Then $S$ and $T$ are
both matrices representing permutations,
each of which consists of a single $n$-cycle. 

\indent 6. Let $C_1$ denote the sequence of integers obtained from
the output of subalgorithm $B$, invoked with input $S$. Similarly
define the sequence $C_2$, obtained from $B$ with input $T$. Then
$C_1$ is of the form $i_1 i_2 \cdots i_n$ and
$C_2$ is of the form $j_1 j_2 \cdots j_n$.

\indent 7. Let $C_3$ denote the reverse of the sequence of integers
$C_1$ (that is, $C_3$ is the inverse of the $n$-cycle $C_1$). Let
$\tau$ denote the output of subalgorithm $D$, invoked with the input
$(C_3,C_2)$. That is, $\tau$ is a sequence of integers
defining  the cyclic decomposition of a permutation, which we also denote $\tau$, 
such that $C_2 = \tau C_3 \tau^{-1}$ (see the output of subalgorithm $B$ for the
format of the integers comprising the cyclic decomposition).  

 \indent 8.  Write to output the sequence of integers $ F(C_1) 0 F(\tau) 0 F(C_3) 0  F(\tau)^R 00 $
 where $F(z)$ denotes the non-$0$ portion of output of subalgorithm $F$, invoked with
 input $z$, and $F(z)^R$ denotes the reverse of the word given by the non-$0$ portion of $F(z)$. Terminate.

\pagebreak

\subsubsection{Subalgorithm $A$}
Convert the sequence of transposions $W$ to its matrix $M$ as given
by the standard matrix representation of $S_n$
defined in \ref{rem:standard_rep}.\\
\\
\indent INPUT: a positive integer $n$, and sequence of transpositions $W$, expressed as a string of
of integers, with the integer $i$ denoting the transposition $(i,i+1)$.\\
\\
\indent OUTPUT: the matrix $M$ representing the permutation given by $W$.\\
\\
\indent ALGORITHM: \\
\\
\indent \indent 1. Let $M$ be the $n \times n$ identity matrix of $GL(n,\Real)$. \\
\indent \indent 2. For each integer $i$ of the input $w$ (reading right to left) \\
\indent \indent \indent 3.  Let $C_1$ be the column of $M$ having $1$ in the $i$-th position \\
\indent \indent \indent 4.  Let $C_2$ be the column of $M$ having $1$ in the $(i+1)$-th position \\
\indent \indent \indent 5.  swap columns $C_1$ and $C_2$ \\
\indent \indent 5. return $M$.\\

\pagebreak

\subsubsection{Subalgorthim $B$} Obtain the cyclic decomposition of $M$. \\
\\
\indent INPUT: a positive integer $n$, and a permutation of $S_n$ given by its matrix $M$ under the standard representation. \\
\\
\indent OUTPUT: a sequence of integers providing the cyclic decompostion of $M$. The end of each cycle
in the decomposition is marked by $0$.  The last cycle in the decomposition is succeeded by
$00$. That is, the output will be of the form: 

%\vspace{-3in}
\begin{displaymath}
i_1 i_2 \cdots i_{k} 0 i_{k+1} i_{k+2} \cdots i_{k+r} 0 i_{k+r+ 1} i_{k+r+2}  \cdots 0 \cdots 0 \cdots i_{n-1} i_n 00
\end{displaymath}
including $1$-cycles, and where the cycles are disjoint.\\
%\vspace{-0.5in}
 \\
\indent ALGORITHM:\\
\\
\indent \indent 1. Let $d$ be an array of $0$s of length $2n$. Let $d(i)$ denote the $i$-th element of the array, with $0$
denoting the first element. Let integers $col, start,i$ all be set to $1$. Let integer $x = 0$. Let $d(x) = i$.\\
\indent \indent 2. Let $x = x+1$. \\
\indent  \indent 3. Let $i$ be the position  (counting from $1$) in which the number $1$ appears in the
column given by $col$.  \\
\indent \indent 4. If $i == start$\\
\indent \indent \indent 5. Let $d(x) = 0$. Let $start$ be the position $j$ (counting from $1$) of the first column
of $M$ such that $j$ does not appear in $d$, or $start = 0$ if no such column.  \\
\indent \indent \indent 6. If $start == 0$ \\
\indent \indent \indent \indent 7. Let $d(x) = d(x+1) = 0$. Return $d$.\\
\indent \indent \indent 7. Else let $col = start$ and $d(x) = col$. Go to $2$. \\
\indent \indent 8. Else let $d(x) = i$. Let $col = i$. Go to $2$. \\

\vspace*{2in}

%\vspace*{-2in}

\pagebreak

\subsubsection{Subalgorithm $C$} Express as a product of two $n$-cycles the even permutation given by the cyclic
decomposition $d$. \\
\\
\indent INPUT:  the parameters $(d,n)$, where $d$ is the cyclic decomposition a permutation of $S_n$.  The expression
$d$ will be of the form:
\begin{displaymath}
i_1 i_2 \cdots i_{k} 0 i_{k+1} i_{k+2} \cdots i_{k+r} 0 i_{k+r+ 1} i_{k+r+2}  \cdots 0 \cdots 0 \cdots i_{n-1} i_n 00
\end{displaymath}
where the $0$s separate the individual cycles and the composition is terminated by $00$ (see the output of subalgorithm
$B$). \\
\\
\indent OUTPUT: two $n$-cycles $c_1,c_2$. Each $n$-cycle will be expressed as a product of (not necessarily disjoint) cycles, where each individual cycle has the form as described in the definition of
the output of subalgoritm $B$. The $n$-cycles $c_1$ and $c_2$ will be separated by the sequence ``$00$''.  The string $c_1 00 c_2$
will be succeeded by $000$. \\
\\
\indent ALGORITHM: \\
\\
\indent \indent 1.  If $d$ consists of exactly one cycle (counting $1$-cycles) of length $n$, that is, if $d$ is of the form
$(i_1,i_2,\cdots i_n)$ where  $n$ is odd \\
\indent \indent \indent /* apply lemma 4 */ \\
\indent \indent \indent 2.  Let $c_1 = c_2 = \underbrace{i_1 i_2 \cdots i_n 0}_{((n+1) / 2) - 1} i_1 i_2 \cdots i_n $ \\

\indent \indent \indent 3. Return $ c_1 00 c_2 000$ \\

\indent \indent 4. Else if $d$ consists of exactly two cycles of lengths $2m$ each, where $n = 4m$, that is, if $d$ is of the form
$i_1,i_3,\cdots, i_{4m-1},0,i_2,i_4,\cdots,i_{4m}$ \\
\indent \indent \indent /* apply lemma 5 */ \\
\indent \indent \indent  5. Let $c_1 = c_2 = i_1 i_2 i_3 i_4 \cdots i_{4m-1} i_{4m}$ \\
\indent \indent \indent 6. Return $c_1 00 c_2 000$. \\

\indent \indent 7.  Else  if $d$ consists of exactly two cycles $\rho_1$ and $\rho_2$ of lengths $2s$ and $2t$ respectively,
with  $s < t$, $n = 2s + 2t$, and $\rho_1 = i_1 i_3 \cdots  i_{4s-1}$ and \\
$\rho_2 =  i_2 i_4 \cdots i_{2s + 2t} i_{4s+1} i_{4s+3} \cdots i_{2s+2t-1}$\\
\indent \indent \indent /* apply lemma 6 */ \\
\indent \indent \indent 8.  Let $c_1 = i_1 i_{4s+1} i_{4s+2} i_{4s+3} \cdots i_{2s+2t} i_2 i_3 \cdots i_{4s}$ \\
\indent \indent \indent 9. Let $c_2 = i_1 i_2 i_3 \cdots i_n$. \\
\indent \indent \indent 10. Return $c_1 00 c_2 000$. \\

\indent \indent 11. Else if $d$ contains a cycle of odd order $d_1$, such that $d_1 = i_1 i_2 \cdots i_k$ \\
\indent \indent \indent /* peel off an odd cycle and apply lemma 4; re-invoke subalgorithm on remainder */\\
\indent \indent \indent 12. Let $d_2$ denote the product of the remaining cycles.  Observe that $d_1$ may be regarded
as an even permutation on a set of $k$ integers $\{i_1,i_2,\cdots i_k\}$, and $d_2$ may be regarded as an even
permutation on a disjoint set of $n-k$ integers.  Let $c_{11} = c_{12} = \underbrace{d_1 0}_{((k+1) / 2) -1} d_1$.  (That is,
we asign to $c_{11}$ and $c_{12}$ the output of $C(d_1,k)$.) \\

\indent \indent \indent 13. set $ (c_{21},c_{22}) = $ the output of subalgorithm $C$, invoked with input \\
$(d_2,n-k)$ .  We have that $c_{21} c_{22}$ is a permutation acting on the set of integers
$\{j_1, j_2, \cdots , j_r \}$, where $\{j_1, j_2, \cdots , j_r \} \bigcap  \{i_1, i_2, \cdots, i_k \} = \emptyset$ \\
and $\{i_1, i_2, \cdots i_k, j_1, j_2, \cdots,  j_r \} = \{1,2,\cdots,n\}$. \\
\indent \indent \indent 14. Let $u$ denote the maximum of  $\{i_1,i_2,\cdots,i_k\}$ and
$v$ denote the maximum of $\{ j_1, j_2, \cdots j_r \}$. \\
\indent \indent \indent 15. Let $c_1 = c_{11} 0  c_{21} 0 u v $. Let $c_2 = u v 0 c_{12} 0 c_{22}$ 
(re-arranging and/or removing $0$, $00$, $000$ delimiters as required). \\
\indent \indent \indent 16. Return $c_1 00 c_2 000$.\\

\indent \indent 17. Else $d$ consists of an even number of cycles of even order \\
denoted $d_1 0 d_2  0 \cdots 0 d_{2k} 00$.\\
\indent \indent \indent /* peel off two even cycles and apply lemmas 5,6.  re-invoke subalgorithm  on the remainder.  */\\
\indent \indent \indent 18. Let $e_1 = d_1 0 d_2 00 $. Let $e_2$ denote the product of the remaining cycles.  Observe that
$e_1$ may be regarded as a permutation on a set of integers $\{i_1,i_2, \cdots , i_{n^\prime} \}$ for $n^\prime < n$. We have that
either $e_1$ is of the form: \\
\indent \indent \indent $e_1 = i_1 i_3 \cdots  i_{4m^\prime-1} 0 i_2 i_4 \cdots i_{4m^\prime}$ \\
where $4m^\prime = n^\prime$, or else $e_1 $ is of the form: \\
\indent \indent \indent $e_1 = i_1 i_3 \cdots  i_{4s^\prime-1} 0  i_2 i_4 \cdots i_{2s^\prime + 2t^\prime} i_{4s^\prime+1} i_{4s^\prime+3} \cdots i_{2s^\prime+2t^\prime-1} 00$ \\
with $2s^\prime + 2t^\prime = n^\prime$.\\
\indent \indent \indent 19. Let $(c_{11}, c_{12})$ be the output of subalgorithm $C$, invoked with \\ input  $(e_1,n^\prime)$. Observe
that with this input, $C$ will return without invoking itself again.\\
\indent \indent \indent 20.  Let $ (c_{21},c_{22}) $ be the output of subalgorithm $C$, invoked with \\ input  $(e_2,n-n^\prime)$. Observe that
$e_2$ may be regarded as a permutation acting on the set of  $n-n^\prime$ integers
$\{ j_1, j_2, \cdots, j_{n-n^\prime}\}$, such that  $\{i_1,i_2, \cdots , i_{n^\prime}\} \bigcap 
\{ j_1, j_2, \cdots, j_{n-n^\prime} \} = \emptyset$ and
 $\{i_1, i_2, \cdots i_n^\prime, j_1, j_2, \cdots j_{n-n^\prime} \} = \{1,2,\cdots,n\}$. \\
\indent \indent \indent  21. Let $u$ denote the maximum of  $\{i_1,i_2,\cdots,i_{n^\prime}\}$ and $v$ denote the maximum of
$\{ j_1, j_2, \cdots, j_{n - n^\prime} \}$.\\
\indent \indent \indent 22. Let $c_1 = c_{11} 0  c_{21} 0 u v $. Let $c_2 =  u v 0 c_{12} 0 c_{22}$ (re-arranging and/ 
or removing $0$, $00$, $000$ delimiters as required). \\
\indent \indent \indent  23. return $c_1 00 c_2 000$.\\
\\
\pagebreak

\subsubsection{Subalgorithm $D$} Conjugator of two $n$-cycles $c_1,c_2$.\\
\\
\indent INPUT: a positive integer $n$, and two $n$-cycles $c_1$ and $c_2$ (which are necessarily conjugate) each expressed in cyclic
decomposition form (see the output of subalgorithm $B$). \\
\\
\indent OUTPUT: the cyclic decomposition (expressed as in the output of subalgorithm $B$) of a 
permutation $\tau$ such that $c_2 = \tau c_1 \tau^{-1}$. \\
\\
\indent ALGORITHM: \\
\\
\indent \indent 1. Let $d$ be an array of $0$s of length $2n$. Let $d(i)$ denote the $i$-th element of the array, 
(counting left to right from $0$). Let $c_i(j)$ denote the value of the $j$-th integer
of the cycle $c_i$ (counting from left to right from $0$). Let integers $val, i$  be set to $0$. Let 
integer $start = c_1(i)$. Let integer $x = 0$. Let $d(x) = start$.\\
\indent \indent 2. Let $x = x+1$ \\
\indent  \indent 3. Let $val = c_2(i)$  \\
\indent \indent 4. If $val == start$\\
\indent \indent \indent 5. Let $d(x) = 0$. Let $start = c_1(i)$, where  $i$ is set to the position (counting left to right from $0$) of the first integer of $c_1$ such that $c_1(i)$ does not appear in $d$, or $start = 0$ if no such column.  \\
\indent \indent \indent 6. If $start == 0$ \\
\indent \indent \indent \indent 7. Let $d(x) = d(x+1) = 0$. Return $d$.\\
\indent \indent \indent 8. Else let $start = c_1(i)$ and $d(x) = start$. Go to $2$. \\
\indent \indent 9. Else let $d(x) = val$. Set $i$ equal to the position (counting left to right from $0$) of the integer of $c_1$ such
that $c_1(i) == val$. Go to $2$. \\

\pagebreak

\subsubsection{Subalgorithm E}

Convert a product of (not necessarily disjoint) cycles to its matrix $M$ as given by the
standard matrix representation of $S_n$ defined in \ref{rem:standard_rep}.\\
\\
\indent INPUT: a positive integer $n$, and a product of (not necessarily disjoint) cycles $S$ in the form
of one member of the pair of permutations as output by subalgorithm $C$.\\
\\
\indent OUTPUT: the matrix $M$ representing the permutation given by $S$.\\
\\
\indent ALGORITHM: \\
\\
\indent 1. Let  $M$ be the $n \times n$ identity matrix of $GL(n,\Real)$. For each cycle $C$ of $S$, reading right to left: \\
\indent \indent 2. Let $i$ be the first integer of $C$ (processing left to right).  
Let $C_p = C_n$ be copies of column $i$ of $M$.  Let $s = i$.\\
\indent \indent 3. If the next integer $j$ of $C$ is  not $0$\\
\indent \indent \indent 4. Let $C_n$ be a copy of column $j$
of $M$. \\
\indent \indent \indent 5. Let column $j$ of $M$ equal column $C_p$. \\
\indent \indent \indent 6. Let $C_p = C_n$.  Let $i$ be the next integer of $C$. Go to 3.\\
\indent \indent  7. Else if the next integer $j$ of $C$ following $i$ is $0$ \\
\indent \indent \indent 8. Let column $s$ of $M$ = $C_p$.  continue with next iteration of the For loop. \\
\indent  9.  Let $T$ be an $n \times n$ matrix of $0$s.  For each integer  $k = 1 \cdots n$ 
\indent \indent 10. Assign to column $k$ of $T$ the number of the column of $M$ in which $1$ appears in position $k$
\indent 11. Let $M = T$. \\

\pagebreak

\subsubsection{Subalgorithm $F$} Express the cyclic decomposition $d$ as a word in the standard generators $\tau_i$
of $S_n$. \\
\\
\indent INPUT:  a positive integer $n$, and the (non-empty) cyclic decomposition $d$ of a permutation of $S_n$.  The expression $d$ will be of the form:
\begin{displaymath}
i_1 i_2 \cdots i_{k} 0 i_{k+1} i_{k+2} \cdots i_{k+r} 0
\cdots i_{k+r+1} i_{k+r+2} \cdots i_n 00
\end{displaymath}
where the $0$s separate the individual cycles and the composition is terminated by $00$ (see the output of subalgorithm
$B$). \\
\\
\indent OUTPUT: a sequence of integers $i_1 \cdots i_j$, terminated with a $0$, where each $i_k$ corresponds to the generator
$\tau_{i_k}$ of $S_n$. The output should be applied right to left. \\
\\
\indent ALGORITHM: \\
\\
\indent 1.  Let integer $i = 0$. Counting from $0$, and reading left to right,
let $d(k)$ denote the  $k$-th integer of the input $d$.  
Let $w$ be an array of $n^3$ integers,  all initialized to $0$.  Let the integer $x = 0$.
We denote by $w(x)$ the  $x$-th entry of $w$,  counting from $0$, and reading left to right.
Let integer $s = d(i)$. Let $i = i + 1$. Go to step $6$.\\
\indent /* check to see if we are done with all cycles */ \\
\indent  2. If $d(i) == 0$ and $d(i+1) == 0$ \\
\indent \indent  3. Return. \\
\indent /* else check to see if we are done with current cycle */ \\
\indent 4.  Else if $d(i) == 0$ \\
\indent \indent /* we are. begin new cycle. $s$ records our starting point */ \\
\indent \indent 5. Let $s = d(i+1)$. Let $i = i + 2$. \\
\indent \indent \indent /* examine the next integer in the cycle. record in $w$ a transposition between
that integer and our our starting point $s$, where the transposition is expressed as a sequence
of adjacent transpositions. */ \\
\indent \indent \indent 6.  Let integer $a = x$. Let integer $c = 0$. Let integer $f$ be the greater of
$s$ and $d(i)$. Let integer $t$  be the smaller of $s$ and $d(i)$. \\
\indent \indent \indent 7. While $ t < f$ \\
\indent \indent \indent \indent 8.   Let $w(x) = t$. Let $x = x+1$. Let $c = c+1$. Let $t = t+1$.\\
\indent \indent \indent /* now, if we had more than one adjacent transposition, repeat all but the last, in inverse order.
the goal: the swap $(m,p)$, where $(p > m)$, is replaced with
$(m,m+1)(m+1,m+2)\cdots(m+p-1,m+p) (m+p-2,m+p-1)\cdots(m,m+1)$ read right to left.  */ \\
\indent  \indent \indent 9.  If $c > 1$ \\
\indent \indent \indent \indent 10.  Let $c = c -1$. While $c >0$.\\
\indent \indent \indent \indent \indent 11. Let $w(x) = w(a + c - 1)$. Let $c = c - 1$. Let $x = x +1$. \\
\indent \indent \indent 12.  Write $w^R$ to output, where $w^R$ denotes the reverse of $w$. Let $i = i + 1$. Clear $w$. Go to step $2$.

\pagebreak

\subsection{Example of application of Cejtin-Rivin algorithm}
We provide an example of the preceding. Let $\sigma \in S_7$ be given by the word 
$\tau_6 \tau_4 \tau_1 \tau_2$. Algorithm $A$ converts this word to a permutation
matrix given by the ordered list $3,1,2,5,4,7,6$, where, counting left to right from $1$,
digit $i$ indicates the position of the integer $1$ in column $i$. Algorithm $B$ converts
this to the cyclic decomposition $(123)(45)(67)$.  We pass $((123)(45)(67),7)$ to
algorithm $C$.  Algorithm $C$ ``peels off''  the odd cycle $(123)$, and invokes itself
twice with inputs $i_1 = ((123),3)$ and (for the remaining cycles) $i_2 =((45)(67),4)$.  To input $i_1$, the
algorithm applies Cejtin-Rivin lemma $4$ to obtain output $c_{11} = c_{12} = (123)(123)$. To input
$i_2$, the algorithm applies Cejtin-Rivin lemma $5$ to obtain output $c_{21} = c_{22} = (4657)$.
(In this example, no further recursion is required. In general, each invocation of algorithm $C$
on the ``remainder'' left after peeling off of one or two cycles, to which Cejtin-Rivin  lemmas $4,5,6$ can be
directly applied, will lead to further invocations.)  Algorithm $C$ then applies Cejtin-Rivin lemma $3$,
with $t = 3$, $s = 4$, to obtain $p_1 = c_{11}c_{21} \nu = (123)(123)(4657)(37)$ and 
$p_2 = \nu c_{21}c_{22} = (37)(4657)(123)(123)$, so that $\sigma \sim p_1 p_2$. We apply algorithms
$E$ and then $B$ to convert $p_1$ and $p_2$ to $n$-cycle form:  $p_1 \sim C_1 =  (1346572)$ and 
$p_2 \sim C_2 = (1746532).$ We then apply algorithm $D$ to $(C_1^{-1} = (2756431),C_2)$ to
obtain $\tau = (12)(45)$ such that $C_2 = \tau C_1^{-1} \tau^{-1}$. We finally apply algorithm $F$
to $C_1,C_3, \tau, \tau^{-1}$ to obtain a word:
\begin{displaymath}
 [F(C_1),F(\tau)] = [\tau_2 \tau_3 \tau_4 \tau_5 \tau_6 \tau_5 \tau_4 \tau_5 \tau_4 \tau_1,\tau_4 \tau_1]  \sim \sigma 
\end{displaymath}
\subsection{Complexity of the Cejtin-Rivin algorthm}

Here we analyze the complexity of the Cejtin-Rivin algorithm.

\begin{prop}[complexity of implementation of Cejtin-Rivin]
The implementation of the Cejtin-Rivin algorithm given above for the
SCPP for $S_n$ executes in time $O(n^2)$, for input of fixed length $k$, and
time $O(k)$, for fixed $n$. (Recall that by \emph{time}, we mean the total number 
executions of pseudo-code statements
of types $1-7$).
\end{prop}
	The relation between $k$ and $n$
can perhaps be described as follows: though for any $\sigma \in S_n$,
and for any positive integer $k$, there exists an \emph{unreduced} word 
representative $w$ of $\sigma$ such that $\text{length(w)} > k$, we also
have that any element of $S_n$ may be represented by a word $w^\prime$
of length at most $n(n-1)/2$.

\noindent \emph{proof}: By \ref{def:equiv_comp}, it
suffices to show that, for fixed $k$, the number of invocations of statements
of types $1-7$ grows as $O(n^2)$, and that for fixed $n$, the number of invocations of 
statements of types $1-7$ grows as $O(k)$. \\
\\
Now, by inspection, the main algorithm requires a fixed finite
number of operations of types $1-4$ and $7$. It must read a string
of $k$ symbols, each in the interval $[0,n-1]$, and pass these as a block to subalgorithm $A$. 
That requires $O(k)$ operations, and is not dependent on $n$. 
	The output of subalgorithm $A$ is a matrix representative of an element of $S_n$, (in the 
	standard representation), which we may implement as an array of integers of size $n$, where each element $i$ of the array is 
an integer $j$ recording the position $j$ in which the $1$ appears in column $i$. This array is passed to subalgorithm $B$,
requiring $O(n)$ operations.  
	Subalgorithm $B$  returns a string of integers bounded by $n^2$ 
(the string is a cyclic decomposition of the permutation given by the input, 
where each cycle has length at most $n$, and where
there are at most $n$ cycles). This is passed to subalgorithm $C$, requiring $O(n^2)$
operations. Subalgorithm $C$ returns a string of size $O(n)$ representing a product of two $n$-cycles 
$c_1$ and $c_2$.   
	These are each passed to subalgorithm $E$, requiring $O(n)$ operations. Subalgorithm $E$
returns strings of length $O(n)$, which are passed to subalgorithm $B$, requiring $O(n)$
operations.
	Subalgorithm $B$ returns $n$-element arrays $C_1$ and $C_2$ representing $n$-cycles. From
this, $C_1^{-1}$ is calculated, requiring $O(n)$ operations. Then $C_1^{-1}$ and $C_2$ are
passed to subalgorithm $D$, requiring $O(n)$ operations. Subalgorithm $D$ returns a 
string of length $O(n)$ representing a permutation $\tau$, which is then inverted, requiring $O(n)$ operations.
	Finally $C_1, C_1^{-1}, \tau$, and $\tau^{-1}$ are passed to subalgorithm $F$, requiring
$O(n)$ operations. The output of each call to subalgorithm $F$ will 
 be string of length $O(n^2)$, which is then written to output; this  will require
 $O(n^2)$ operations.
 	Hence the main algorithm requires a  number of operations of types $1-7$
that is at worst $O(n^2)$ and $O(k)$ (not counting the number of operations required by each 
subalgorithm).	Now we must verify that the same can be said of each subalgorithm (each of
which is invoked a fixed number of times by the main algorithm). \\

Subalgorithm $A$ accepts an input $w$ consisting of pairs of integers
denoting transpositions. For each operation $w$, it must interchange
two $n$-digit strings of symols. Therefore; for a fixed $n$, it is linear in $w$.\\

Subalgorithm $B$ involves inspection of each cell of the $n \times n$ array that it
accepts as input, coupled with the writing of a string of length bounded by $2n$. Its
execution time is therefore  $O(n^2)$.\\
 
Subalgorithm $C$ will require a number of operations of types $1-5$
that is linear in $n$, if the cyclic decomposition of the
permutation it accepts as input consists of exactly one or two cycles.
Otherwise, the subalgorithm invokes
itself, after which it performs a number of operations of
types $1-5$ that is linear in $n$.  Each invocation of $C$ requires
a number of operations that is linear in $n$.  The number of self-invocations required grows linearly with $n$.  
Hence the subalgorithm is  $O(n^2)$. It does not depend on $k$.\\

Subalgorithm $D$ must compare two $n$-digit strings that it accepts as input.  It
must separately perform $n$ searches of a $3n+1$-digit array. It is therefore
quadratic  in $n$. \\

Subalgorithm $E$ is the inverse of subalgorithm $B$, and involves the reading
of a string of length $j$ that is a linear multiple of $n$, coupled with, 
for each symbol of $j$, the exchange of a fixed number of
integers in a list of length $n$.  It then scans the list $n$ times to create a new list of length $n$,
which is then copied to the original list. Its execution time is then $O(n^2)$. It does not depend
on $k$. \\

Subalgorithm $F$ accepts an input bounded by $Ln$ for some constant
$L$. It requires $n$ iterations of a sequence of operations of types $1-5$,
where number of steps in the sequence is also bounded by $Dn$ for
some constant $D$. Hence the algorithm is quadratic in $n$. \\

Since each subalgorithm is at worst linear in $k$ and at worst quadratic in $n$, 
by the argument in the beginning of this proof, the algorithm
itself is, too. Note that subalgorithms $B$ and $D$ could be improved by representing 
the elements of the standard representation of $S_n$,
not by matrices, but by arrays of integers, such that array entry
$i$ contains the number $j$ if the column $i$ of the corresponding matrix
contains a $1$ in the $j$-th position. (It is less clear
how to linearize subalgorithm $C$, as it must recurse $O(n)$
times, invoking $O(n)$ operations for each recursion.) \indent $\Box$

\pagebreak

\subsection{Word reduction in $S_n$}\label{subsec:symm_reduce}

We shall have need in what follows for a technique by which one can reduce an arbitray word $w$ representing an
element $\sigma \in S_n$ to a unique word $\rho(\sigma)$ that will we call the \emph{canonical representative} of
$\sigma$.

\begin{defn}[\textbf{ShortLex} order]
Let the alphabet $\Sigma$ in definition \ref{def:alphabet}  be provided with a \emph{total order}; that is, a 
reflexive, transitive, antisymmetric relation $ < $ such that for all $a,b \in \Sigma$, either $a < b$ or
$b < a$.  The we can define the $\textbf{ShortLex}$ order $<$ on $\Sigma^*$, where for any two strings
$v,w \in \Sigma^*$, we have $v < w$ if an only if $v$ is shorter than $w$, or else $v$ and $w$ are of the
same length, and $v(i) < w(i)$, where the $v(i)$ and $w(i)$, the $i$-th letters of $v$ and $w$, are the first at 
which the two string differ (reading right to left and counting from one). The \textbf{ShortLex} order is
a well-ordering; any subset of $\Sigma^*$ will have a unique smallest element given by the 
\textbf{ShortLex} order. \cite{Epstein}
\end{defn}

\begin{ex}
Let $\Sigma_n$ be the alphabet $\Sigma_n = \{1,2,3, \cdots,  n -1\}$. Then we can impose the total 
order $1 < 2 < 3 < \cdots < n - 1$ on $\Sigma_n$. Then for strings $13243, 13245 \in \Sigma_6^*$, we
have that $13243 < 13245$.
\end{ex}

\begin{rem}[the group $S_n$ as a language]\label{rem:symm_as_language}
Associate to the symbol $i \in \Sigma_n$ the generator $\tau_i$ of $S_n$.  We can then regard 
the group $S_n$ as the group of equivalence classes of strings over $\Sigma_n$, where two strings
$u$ and $v$ are equivalent if and only if we can transform the string $u$ to the string $v$
by repeated application of the following identifications: 
\begin{equation}\label{eq:symm_id}
(i i \sim e), ( (i+1) \, i \, (i+1) \sim i \, (i+1) \, i), ((i+k) \, i \sim i \, (i+k) \{ \text{for} \,\, k > 1 \} ).
\end{equation}
where $e$ denotes the empty string. (That is, $S_n$ is the group having generators  $1,2,\cdots n-1$ and
relations $i i = e, (i+1)\,i\,(i+1) = i \, (i+1)\, i, (i+k) \, i = i \, (i+k) \{\text{for}\,\, k > 1 \}$). 
\end{rem}

\begin{defn}[reducible]
Observe that the identifications $(u \sim v)$ given in \ref{eq:symm_id} are listed with the \textbf{ShortLex}
smaller term on the right.  We say that a string $u \in \Sigma_n^*$ is \emph{reducible} to a string $v \in \Sigma_n^*$,
denoted $u \rightarrow^* v$, if by repeated application of the identifications (\ref{eq:symm_id}) $u$ can be 
transformed to $v$, and such that $v < u$ in the \textbf{ShortLex} order.
\end{defn}

Given a set of rules $S$ such as the above, but such that each pair of rules $(u_1 \sim v_1),(u_2 \sim v_2)$ satisfy the
following  two properties (where $e$ denotes the empy string):
\begin{equation}
\begin{array}{l}
\text{if}\, u_1 = rs \,\, \text{and} \,\, u_2 = st \,\, \text{with} \,\,s,t,r \in \Sigma_n^{*}, \,\text{and}\, s \neq e,\\ 
\text{then there exists} \,\,  w \in \Sigma_n^{*} \,\, \text{such that}\,\,  v_1 t  \rightarrow^{*} w \,\, \text{and}\,\, r v_2 \rightarrow^* w
\end{array}
\end{equation} 

\begin{equation}
\begin{array}{l}
\text{if}\, u_1 = rst \,\, \text{and} \,\, u_2 = s \,\, \text{with} \,\,s,t,r \in A_n^{*}, \,\text{and} \, s \neq e \\ 
\text{then there exists} \,\,  w \in \Sigma_n^{*} \,\, \text{such that}\,\,  u_1   \rightarrow^{*} w \,\, \text{and}\,\, r v_2 t \rightarrow^* w
\end{array}
\end{equation} 

\noindent  one can reduce any word $w$ representing an element of $S_n$ to a canonical form $r(w)$, such that for each $\sigma \in S_n$, 
we have a unique representative $r(w)$, obtainable from any word $w$ representing $\sigma$. This reduction is given
by the following algorithm \cite{Epstein}, pp.116-126:\\
 
\subsubsection{Reduction algorithm $\mathcal{A}$ for $S_n$\\}\label{subsub:reduce} 

\medskip

\indent \emph{PROBLEM: Reduction of $w$ to its canonical form.} \\
\indent INPUT: a word $w \in \Sigma_n^*$ \\
\indent OUTPUT: $\mathcal{A}(w)$ \\
\indent ALGORITHM: \\

\indent 1. Select the subset of rules $S_1 \subset S$ such that for each rule, the left-hand side matches a substring of $w$. \\
\indent 2. If $S_1 = \{ \emptyset \}$, terminate.\\
\indent 3. Let $S_2 \subset S_1$ be the set of rules in $S_1$, such that for each rule, the first character of the
matching substring of $w$ is as far to the left of $w$ as possible.\\
\indent 4. Let $S_3 \subset S_2$ be the set of rules in $S_2$ have the smallest left-hand sides, in \textbf{ShortLex} order.\\
\indent 5. Let $r \in S_3$ be the rule having the shortest right-hand side, in \textbf{ShortLex} order. \\
\indent 6. Apply rule $r$ to the left-most possible substring of $w$. Go to step $1$.\\
\\
Such a set of rules $S$ constitute what is known as a  \emph{complete set $S$ of Knuth-Bendix rules}.  The relations
given above for $S_n$ (\ref{eq:symm_id}) are \emph{not} a complete set of Knuth-Bendix rules. However, given a finitely presented group $G$ having generators $X = \{x_1, x_2, \cdots, x_n \}$ and relations 
$R= \{(u_1,v_1), (u_2, v_2), \cdots (u_m, v_m) \}$, such that $u_i > v_i$ with respect to a \textbf{ShortLex} order, one can attempt to obtain a complete set of Knuth-Bendix rules as follows \cite{Epstein}, pp.116-126. The procedure will always terminate for a finite group $G$ \cite{Holt}: \\

\indent 1. start with the set of rules $S_0 = S$ given by the $(u_i, v_i)$. \\
\indent 2. For each pair of rules $r_1$, $r_2 \in S_i$ such that the left-hand sides  overlap (are either of the form  
$r_1 = (uv,x)$ and $r_2 = (vw,y)$ with $v \neq e$, or $r_1 = (v,x)$ and $r_2 = (uvw,y)$: \\
\indent \indent 3. If no such overlap, terminate. \\
\indent \indent 4. Else perform the reduction process desribed above on the word $uvw$,  once starting with rule $r_1$,
to obtain a word $t_1$, and once starting with $r_2$, to obtain a word $t_2$.\\
\indent 5. If $t_1 \neq t_2$ (say $t_1 < t_2$), add the rule $(t_2,t_1)$ to $S_i$ to obtain a  new set  $S_{i+1}$.\\
\indent 6. Go to 2.

%\vspace{-0.5in}

\begin{rem}[complexity of the reduction process]
In what follows, we assume the existence of efficient implementations of algorithms for word reduction in $S_n$, (that is, polynomial time for fixed $k$ and varying $n$, and also polynomial time for fixed $n$ and varying $k$), that are functionally equivalent to the preceding (see for example the programs listed in \cite{Holt} p. 8). 
\end{rem}

\pagebreak

\section{The Braid Groups $B_n$}

In this section we present information about the braid groups $B_n$, with a focus on braids as equivalence
classes of words in the standard Artin generators and their relations.

\subsection{The braid groups $B_n$: algebraic, geometric, and topological definitions}

\begin{defn}[Braid Group: algebraic definition, \cite{Ohtsuki},  p.23] \label{def:braid-algebraic}
The braid group on n strands (denoted $B_n$) is given by the
generators $\sigma_1 \cdots \sigma_{n-1}$ subject to the relations
\begin{equation}
\nonumber
\begin{array}{lr}
\sigma_i \sigma_j = \sigma_j \sigma_i & \text{if} \mid i - j \mid > 2,\,\, 1 \leq i,j \leq n - 1 \\
\sigma_i \sigma_{i+1} \sigma_i = \sigma_{i + 1} \sigma_{i}
\sigma_{i + 1} & 1 \leq i \leq n - 2.
\end{array}
\end{equation}
\end{defn}

\begin{defn}[Braid Group: geometric definition, \cite{Ohtsuki}, p.22] \label{def:braid-geometric}
A braid on $n$ strands may be regarded geometrically as the union of
$n$ arcs embedded in $\Real^2 \times [0, 1]$, such that the boundary
is the set $\{1,2,\cdots,n\} \times \{0 \} \times \{0,1\} \in
\Real^2 \times [0,1]$ and such that no arc has a critical point with
respect to the boundary (see figure $1$). Two braids are isotopic if
they are related by an isotopy of $\Real^2 \times [0, 1]$ preserving
the boundary and the vertical coordinate.  Any geometric braid is
obtained by attaching vertically a series of elementary braids, with
each elementary braid  given by an elementary braid diagram (which is the projection
of the braid onto the $x-z$ plane, w/the $y$-axis directed vertically outwards
from the face of the page). Each elementary braid corresponds
to a generator of the braid group (see figure $2$).  Isotopy introduces an equivalence 
relation on the set of braids, and the set of equivalency classes is a group, with action given
by the attaching process just described. This group is in fact
isomorphic to the group of equivalency classes of words in the
$\sigma_i^{\pm}$ given in the preceding definition see
(\cite{Murasugi}). Finally, a  $\emph{pure braid}$ is one for which each
each strand $i$ beginning at $(x_i,0) \in \Real^2 \times [0,1]$ ends
at $(x_i,1) \in \Real^2 \times [0,1]$.
\end{defn}

\begin{defn}[Braid Group, topological definition, \cite{Kassel2}, p.8]
A third definition of the braid group, equivalent to the preceding two, may
be given as follows.  Let $Y_n$ denote the set of all $n$-tuples
$ (z_1, \cdots ,z_n)$ of points in $\Complex$ such that $z_i \neq z_j$
for $i \neq j$.  Let $X_n = Y_n / S_n$: this is the \emph{configuration
space} of $n$ (unordered) points in $\Complex$. Select a base point
$ p \in X_n $. Then $B_n$ is the fundamental group $\pi_1(X_n,p)$ of
$X_n$.
\end{defn}

\begin{defn}[Braid word, word length, \cite{Murasugi}, pp.15-19]
Let $b$ be a braid represented by a word $w$ in the braid group
generators $\sigma_i$ and their inverses.  We will refer to such $w$
as \emph{braid words}.  The number of $\sigma_{i_j}^{\pm}$ appearing
in the word $w$ is said to be the \emph{word length} of $w$, denoted
$l(w)$.  If $w$ is a word of miminal length representing $b$, then
$\mid b \mid = l(w)$  is said to be the \emph{length} of $b$ (observe
that the braid relations preserve word length, so any word of minimal
length representing $b$ determines the length of $b$).

\end{defn}
For example, if $b$ is the braid word $\sigma_1 \sigma_1 \sigma_2$, then
$\mid b \mid = l(w) = 3$.

\subsection{A homomorphism from $B_n$ to $S_n$}

\begin{rem}[projection of $B_n \rightarrow S_n$, \cite{Epstein}, ch. 9]
A braid $b$ maps the set of points $\{1,2,\cdots,n\} \times \{0\} \times \{0\}$ to the set of
points $\{1,2,\cdots,n\} \times \{0\} \times \{1\}$. The map is $1-1$, and
the point $\{i\} \times \{0\} \times \{0\}$ is taken to the point $\{b(i)\} \times \{0\} \times \{1\}$.
Hence, a braid defines a permutation:
\begin{equation}
\begin{array}{cccc}
1 & 2 & \cdots & n \\
b(1) & b(2) & \cdots & b(n)
\end{array}
\end{equation} \\
\\
In other words, letting $\tau_i$, for $i = 1, \cdots, n-1$, denote the standard generators of the
symmetric group $S_n$ (where $\tau_i$ denotes the permutation of the $i$-th and $i+1$-th elements),
one obtains a surjective group homomorphism
$\rho: B_n \rightarrow S_n$ given by $\rho(\sigma_i) \rightarrow \tau_i$.  
\end{rem}

\subsection{Right greedy normal form}

The information in this section is taken from chapter 9 of \cite{Epstein}. 

\begin{defn}[permutation braids, \cite{Epstein}, pp.182-190]

A \emph{positive} braid is a braid consiting only of positive crossings (that is, crossings represented
only by the $\sigma_i$, and not by any $\sigma_i^{-1}$). A \emph{permutation braid} is a positive braid 
in which no two strands cross twice.  We denote the permutation braids of $B_n$ by $D_n$.
\end{defn}

\begin{rem}[isomorphism of $D_n$ and $S_n$, \cite{Epstein}, pp.182-190] \label{rem:dn_iso_sn}
It is known that the permutation braids $D_n$ are in $1-1$
correspondence with the elements of the symmetric group $S_n$. The
isomorphism is given by the restriction to $D_n$ of the map $\rho:
B_n \rightarrow S_n$  defined above.
\end{rem}

\begin{rem}[permutations given by sets of labels exchanged under permutation,
\cite{Epstein}, pp.182-190]\label{rem:permord}
Consider the action of an element $\sigma$ of $S_n$ on the set of integers $\{1,2,\cdots,n\}$.  Let
 let $R_\sigma = \{(i,j) | i < j\, \text{and}\, \sigma(j) < \sigma(i)\}$.  This set characterizes the
 permutation $\sigma$. We can define a partial order relation $\geq$ on $S_n$ by letting
$\sigma \geq \tau$ if $R_\sigma \supseteq R_\tau$
 (recall a partial order relation on a set $S$ is a binary operation on $S$ that is reflexive,
 antisymmetric, and transitive). \\
 \\
 Given two permutations $\sigma$ and $\tau$, and the order relation $\geq$, there is a unique maximal element
$\sigma \wedge \tau$ such that $\sigma \geq \sigma \wedge \tau$ and
 $\tau \geq \sigma \wedge \tau$. The element is given by the
 recursive formula:
\begin{equation}
 R_{\sigma \wedge \tau} = \{ (i,k) \in  R_{\sigma} \cap R_{\tau} | (i,j) \in R_{\sigma \wedge \tau}\,\, \text{or}\,\,(j,k) \in R_{\sigma \wedge \tau}\, \forall \, j\, \text{s.t.}\, i < j < k\}.
\end{equation}\\
\\
(Observe that for starters, any pairs $(i,i+1) \in  R_{\sigma} \cap
R_{\tau} $ will also appear in $R_{\sigma \wedge \tau}$). It turns
out that there is also a unique minimal element $\sigma \vee \tau$
larger than both $\sigma$ and $\tau$).
\end{rem}

\begin{defn}[heads and tails]
Let $P_n$ be the set of positive braids of $B_n$. Let $a,b,c \in P$. If $ab = c$, we say that $a$ is a \emph{head}
of $c$ and $b$ is a $\emph{tail}$ of $c$ (one writes $a \prec c$ and $c \succ b$).  It turns out (see \cite{Epstein}
p. 188) that for braids $b,c \in D_n$, we have  $c \succ b$ iff $\rho(c) \geq \rho(b)$.
\end{defn}

\begin{rem} [facts about $D_n$, \cite{Epstein}, pp.186-187] \label{rem:permfacts}
 There is a unique maximal element of $D_n$ -- the element $m \in D_n$ such that $R_{\rho(m)} \supseteq R_{\rho(l)}\,\forall\,l \in D_n$
. It is known as the ``Garside element,'' denoted $\Omega$. It is
given by $\Omega =
 ( \sigma_{1} \sigma_{2}\cdots\sigma_{n-1})(\sigma_{1} \sigma_{2} \cdots \sigma_{n-2})\cdots(\sigma_1)$,
 which is a  $180$-degree twist of all strands. Here and in the remainder of this section, we apply braid words from
right to left (treating the word as an operator on $B_n$), and draw them horizontally, moving right to left,
with strand $1$ on the top; for positive crossings, the $i$-th
strand will pass over the $i+1$-th. Hence the Garside element of
$D_3$ is $(\sigma_1 \sigma_2) \sigma_1$ (see figure \ref{fig:3}). We
have the relations $\Omega \sigma_{i} = \sigma_{n-i} \Omega$, and
$b\Omega^2 = \Omega^2 b$ for all $b \in B_n$.  We choose for each
element $d_k$ of $D_n$ a representative braid word $w_k$ and call
this the \emph{canonical representative} of $d_k$.  For instance,
for $D_3$ we may choose the representatives $ \{ e,\sigma_1,
\sigma_2, \sigma_1 \sigma_2, \sigma_2 \sigma_1, \sigma_1 \sigma_2
\sigma_2 \} $, where $e$ denotes the braid word consisting of no
letters.  In what follows, we shall always choose canonical representatives
corresponding to the canonical representatives of elements of $S_n$ as
described in subsection \ref{subsec:symm_reduce}.

\end{rem}

%\vspace{-0.205in}

\begin{defn}[right greedy normal form (r.g.n.f.), \cite{Epstein} pp.191-196]\label{def:rgnf}
A braid word $w$ is in \emph{right greedy normal form} if it is  written
in the form: \\
\begin{equation}\label{eq:rgnf1}
w = w_1 w_2 \cdots w_m \Omega^{p}
\end{equation}\\
\\
where each $w_k$ is the canonical representative of an element of
$D_n$, such that if two strands that are adjacent at the boundary of
$w_k$ and $w_{k-1}$ cross in $w_{k-1}$, then they cross in $w_k$,
and where $p$ may be positive or negative, and where none of the
$w_k = \Omega$.  The integer $p$ in (\ref{eq:rgnf1}) is called the
\emph{infimum} of $w$, and the integer $m$ is called the
\emph{canonical length} of $w$. In what follows, we denote the right greedy
normal form of a braid word $a$ by $r.g.n.f.(a)$. It is known that a classical
Turing machine exists to convert a braid word $w$ to $r.g.n.f.(w)$ in time
quadratic in the length of $w$.
\end{defn}

%\vspace{-0.205in}

\begin{rem}[Existence and Uniqueness of r.g.n.f. \cite{Epstein}, p.195]\label{rem:rgnf_unique}
There is a unique braid word in r.g.n.f. for each braid $b \in B_n$.
\end{rem}

\begin{rem}[Method of writing a braid in right-greedy normal form \cite{BirmanBrendle}]\label{rem:rgnf2}
Given a braid word $w$, having $k$ letters $\sigma_i^{-1}$, we may transform it into r.g.n.f as follows: \\

1.  For each $\sigma_i^{-1} \in w$, multiply $w$ on the right by $\Omega^2 \Omega^{-2}$.  Call this new braid word
$w^\prime$.\\
\\
2. Using the commutivity of $\Omega^2$, move a copy of $\Omega^2$ to the immediate left
of each $\sigma_i^{-1}$. This replaces each $\sigma_i$ with a positive word of the form $\Omega(\bullet)$, where
$(\bullet)$ consists of $\mid \Omega \mid - 1$  letters. Call this braid word $w^{\prime \prime}$.\\
\\
3. Using \ref{rem:permfacts}, move all copies of $\Omega \in w^{\prime \prime}$ to the immediate left of the
$\Omega^{-2k}$ term added in step $1$. This yields a word of the form:\\
\\
\begin{equation}
w^{\prime \prime \prime} = \sigma_{p_1} \sigma_{p_l} \cdots \sigma_{p_m} \Omega^r
\end{equation}
where $r$ is positive or negative, and all the $\sigma_{p_j}$ are positive twists.  \\
\\
4. Next we want to express the subword $p = \sigma_{p_1} \cdots
\sigma_{p_n}$ as product of canonical representatives of $D_n$. We
proceed as follows. Scanning from right to left, we find the largest
substring $s$ of letters of $p$, $s = \sigma_{p_i} \cdots
\sigma_{p_m}$, such that $s$ represents a tail of the braid
represented by $p$, and such that no two strands of $s$ cross twice.
Replace the substring with its canonical representative $w_k \in
D_n$. Repeat this procedure iteratively beginning with the remaining
$\sigma_{p_1} \cdots \sigma_{p_{i-1}}$ letters, to obtain a braid
word:

\begin{equation}\label{eq:rgnf2}
w^{(4)} = w_1 w_2 \cdots w_q \Omega^r
\end{equation}
\\
5. We have that each $w_k$ represents a braid that is a tail of the
braid represented by the $w_1 \cdots w_{k-1}$. We want the $w_k$ to
have the property that if two strands that are adjacent at the
boundary of $w_k$ and $w_{k-1}$ cross in $w_{k-1}$, then they cross
in $w_k$. It turns out that this is equivalent to saying (see
\cite{Epstein} pp.190-196) that each $w_k$ is \emph{maximal} tail
for the word $w_1 \cdots w_{k}$. This means that $w_k$ represents
the maximal element of $D_n$, such that the element is a tail of the
braid given by $w_1 \cdots w_{k}$. (That is, if $v_j$ is any other 
element of $D_n$ that is a tail of $w_1 \cdots w_k$,  then $\rho(w_k) \geq \rho(v_j)$).
(We can define a ``minimal head'' of $w_1 \cdots w_{k}$ analogously). \\
\\
We can obtain a word $w^{(5)}$ in such a form as follows. Consider
the product $ab$ of two canonical representatives $a$ and $b$ of
elements of $D_n$. Let $b^R$ denote the reversal of $b$, that is, if
$b = \sigma_{i_1} \cdots \sigma_{i_l}$, then $b^R = \sigma_{i_l}
\cdots \sigma_{i_1}$.  Let $m(a,b)$ be the canonical representative
of the element mapped by $\rho$ to the element $\rho(a) \wedge
(\rho(\Omega) \rho(b^R))$ of $S_n$. It turns out (see
\cite{Epstein},
p. 191, Prop. 9.2.1) that $m(a,b)$ is the unique maximal tail of $a$, such that $m(a,b)b \in D_n$.  \\
\\
So, we proceed as follows. We consider the product $w_1 w_2$, let
$t_0 = w_1$, let $w_{q+1} = e$, and let $h_1 = t_0 
m(t_0,w_2)^{-1}$. That is, $h_1$ is the minimal head of $t_0 w_2$.
Let $t_1 = m(t_0,w_2)w_2$. Next we let $h_2 = t_1  m(t_1,w_3)^{-1}$.
Then let $t_2 = m(t_1,w_3)w_3$. In general we have
\begin{equation}
h_i = t_{i-1}m(t_{i-1},w_{i+1})^{-1}\,,\,\,t_i =
m(t_{i-1},w_{i+1})w_{i+1}
\end{equation}
We obtain a sequence
\begin{equation}\label{eq:rgnf3}
w^{(5)} = h_1 h_2 \cdots h_{q-1} t_q \Omega^r
\end{equation}
where each $h_i$  satisfies the properties of \ref{def:rgnf}, as
does $t_q$.
\end{rem}
\
\begin{ex}[computing r.g.n.f. for a simple commutator of elements of $D_3$]
An example is in order. The elements of $D_3$ have represenatives as
follows:\\
 $e,\sigma_1, \sigma_2,\sigma_1 \sigma_2, \sigma_2 \sigma_1,
\sigma_1 \sigma_2 \sigma_1$.
Under $\rho$, these elements map to permutations given by the
following $R$-sets:\\
\begin{equation}
\begin{array}{lc}
R_{\rho(\sigma_1)} = &\{(1,2) \} \\
R_{\rho(\sigma_2)} = &\{(2,3) \} \\
R_{\rho(\sigma_1 \sigma_2)} = & \{ ((1,2),(1,3) \} \\
R_{\rho(\sigma_2 \sigma_1)} = & \{ ((1,3),(2,3) \} \\
R_{\rho(\Omega)} = & \{ (1,2),(1,3),(2,3) \} \\
\end{array}
\end{equation}
The complimentary permutations $\rho(\Omega)\rho(\sigma)$ are given by:
\begin{equation}
\begin{array}{lc}
R_{\rho(\Omega) \rho(\sigma_1)} = &\{(1,3),(2,3) \} \\
R_{\rho(\Omega) \rho(\sigma_2)} = &\{(1,2),(1,3) \} \\
R_{\rho(\Omega) \rho(\sigma_1 \sigma_2)} = & \{ (1,2) \} \\
R_{\rho(\Omega) \rho(\sigma_2 \sigma_1)} = & \{ (2,3) \} \\
R_{\rho(\Omega)\rho(\Omega)} = & \emptyset \\
\end{array}
\end{equation}
The reversals of each element are given by $\sigma_i^R = \sigma_i, (\sigma_1 \sigma_2)^R =  \sigma_2 \sigma_1, \Omega^R = \Omega$.\\
\\
So, say we are given a braid word $w = [\sigma_1,\sigma_2] = \sigma_1 \sigma_2 \sigma_1^{-1} \sigma_2^{-1}$. We apply steps $1$
and $2$ to obtain:
\begin{equation}
 w^{(2)} = \sigma_1 \sigma_2 \Omega^2 \sigma_1^{-1} \Omega^2 \sigma_2^{-1} \Omega^{-4}  = \sigma_1 \sigma_2 (\Omega \sigma_1 \sigma_2) (\Omega \sigma_2 \sigma_1) \Omega^{-4}.
\end{equation}
Now we move the $\Omega$'s to the right to obtain $w^{(3)} =  \sigma_1 \sigma_2 \sigma_2 \Omega^{-1}$. Working right
to left, we then factor this as $w^{(4)} = (\sigma_1 \sigma_2) (\sigma_2) \Omega^{-1}$.  Then we 
set  $t_0  = (\sigma_1 \sigma_2)$. Now 
$m(\sigma_1 \sigma_2,\sigma_2)  = \rho(\sigma_1 \sigma_2) \wedge (\rho(\Omega) \rho(\sigma_2^R)) $. This
is the wedge of the sets
$\{ (1,3),(2,3) \}$ and $\{ (1,2),(1,3) \}$, which is $ \emptyset$, and this corresponds to the canonical representative $e$.
So $h_1 = (\sigma_1, \sigma_2) e^{-1} = (\sigma_1 \sigma_2)$ and $t_1 = e \sigma_2 = \sigma_2$. Then $m(t_1,w_3) = 
m(\sigma_2,e) = \sigma_2$ so $h_2 = t_1 m(\sigma_2,e)^{-1} = e$ and $t_2 = m(\sigma_2,e) w_3 = \sigma_2$ So
we have $w^{(5)} = h_1 t_2 = (\sigma_1 \sigma_2) \sigma_2$. (In this particular case, we were not able to pull any additional crossings to the right in step $5$.)
\end{ex}

\begin{rem}[right greedy normal form is in $P$]\label{rem:rgnf_in_P}
It is known (\cite{Epstein}, \cite{BirmanBrendle}) that an algorithm exists to convert
an arbitrary braid word to right greedy normal form in time quadratic in the
length of the braid word. From this together with (\ref{rem:rgnf_unique}),
it follows that there is a polynomial time algorithm for solving the word
problem in $B_n$; that is, there is a classical Turing machine that
can compare two arbitrary braid words $w_1, w_2$, for fixed $n$, in time polynomial in the
sum of the lengths of $w_1$ and $w_2$. 
\end{rem}

\begin{rem}[positive cancellation]\label{rem:cancellation}
Positive braids obey a right cancellation law: for positive braids
$a,a^\prime,b$, we have that $a^\prime b = a b$ iff $a = a^\prime$. A similar
law holds for left cancellation.  Note that in what follows, we will
write, for braid words $w_1, w_2$, that $w_1 = w_2$ when
$w_1$ and $w_2$ are identical as braid words (that is, they consist
of exactly the same sequence of $\sigma_i^{\pm}$, and $w_1 \sim w_2$
when $w_1$ and $w_2$ representative equivalent braids, but are not
necessarily identical as braid words (that is, they are equivalent as braid
words under the braid relations, but may have different ``spellings''.
For example, we will write $w_1 \sim w_2$ if $w = \sigma_1 \sigma_2 \sigma_1$
and $w_2 = \sigma_2 \sigma_1 \sigma_2$. We will write $w_1 = w_2$ if
both are the strings $\sigma_1 \sigma_2 \sigma_1$. Given a braid word $w$,
$\bar{w}$ will denote the corresponding braid.

\end{rem}

\pagebreak

\section{The SCPP for $B_n$: Extension of Cejtin-Rivin}

In this section we describe an extension of the Cejtin-Rivin algorithm to a subset of the
permutation braids of $B_n$.

\subsection{The simple commutator promise problem for $B_n$}

\begin{defn}[simple commutator promise problem]
The \emph{simple commutator promise problem (SCPP)} for the braid
groups is defined as follows: Let  $w$ be a braid word representing a braid
$\gamma \in B_n$, such that $\gamma$ is of the form $\gamma = [a,b] := aba^{-1}b^{-1}$
for braids $a$ and $b$.  Find braid words $x$ and $y$ such that
$[x,y]$ represents the braid $\gamma$; that is, $[x,y]  \sim w$.
\end{defn}

\begin{defn}[related conjugacy problem]
One active area of braid group research is the \emph{conjugacy
problem}: given two braid words $w$ and $w^\prime$, can one
determine whether a braid word $x$ exists such that $w \sim x
w^\prime x^{-1}$? That is, do the two braid words represent conjugate braids?
No polynomial time solution is currently known for this problem, though
exponential time algorithms exist (see \cite{BirmanBrendle}).
Related to this is the ``search'' version of the conjugacy problem:
given braid words $g$ and $h$ that represent conjugate braids, can
one find a braid word $x$ such that $g \sim x h x^{-1}$? As with the
conjugacy problem itself,  no known
polynomial time algorithm for the conjugacy search problem (CSP) currently
exists.
\end{defn}

\begin{rem}\label{rem:scpp_2_csp}
Naturally, the simple commutator promise problem and the conjugacy
search problems are related. One way to see this: say $w$ is a braid
word representing a simple commutator $\gamma$. Say one had a method
of finding a braid word $a$ such that $w \sim  [a,x]$ for some
braid word $x$. Then to solve the SCPP, one would need to find a braid
word $x$ satisfying:

\begin{equation}\label{eq:csp}
x a^{-1} x^{-1} \sim a^{-1} w
\end{equation}.
\\

\end{rem}

\subsection{An extension of the Cejtin-Rivin algorithm to the SCPP for permutation braids}

We would like to extend the Cejtin-Rivin algorithm described above to
solve the SCPP for braids that are simple commutators of
permutations braids; that is, for all braids of the form $[d_1,d_2]$
for $d_1,d_2 \in D$. We will denote this set of braids by
$[D_n,D_n]$, though we must remember that the set does not have the
structure of a subgroup of $B_n$, since the product of two elements
of $[D_n,D_n]$ need not be in $[D_n,D_n]$. \\

\begin{defn}[inverse of $\rho:B_n \rightarrow S_n$]\label{def:rho_inv}

As we saw above, we have that the set of permutation braids $D_n$ is
in $1-1$ correspondence with the elements of the group $S_n$, with
the correspondence given by the natural projective homomorphism
 $\rho: B_n \rightarrow S_n$, when restricted to $D_n$. 
 For $\sigma \in S_n$, let $\rho^{-1}(\sigma)$ denote the unqiue element
  $d \in D_n$ such that $\rho(d) = \sigma$. \\
  
\end{defn}

Now, let $b$ be an element of the set $[D_n,D_n]$. Then 
$b = [d_1,d_2]$ for $d_1, d_2 \in D_n$. Since $b$ is a simple commutator,
any expression of $b$ as a word in the standard Artin generators of
$B_n$ will contain an even number of generators, as the braid
relations preserve word length. Hence $\rho(b) \in A_n$, the
commutator subgroup of $S_n$.  (Indeed, by the same argument,
for \emph{any} element $b$ of the commutator subgroup $[B_n,B_n]$,
we have that $\rho(b) \in A_n$.)   Let $CR(\sigma)$ be the result of
applying the Cejtin-Rivin algorithm described above to the
permutation $\sigma \in A_n$. So we have that $CR(\rho(b)) = [\psi_1,
\psi_2]$ for permutations $\psi_1, \psi_2$ of $S_n$. We would like
to know: when is it the case that $b = [\rho^{-1}(\psi_1),\rho^{-1}(\psi_2)]$? 
We denote the subset of all such  elements in $[D_n,D_n]$ by $K$.
\\

\begin{defn}[the set of braid words $\mathcal{K}$]\label{def:kwords}
For each $\sigma \in A_n \subset S_n$, let $CR(\sigma)$ be the 
simple commutator $[\psi_1(\sigma),\psi_2(\sigma)]$ that is the output of the
Cejtin-Rivin algorithm applied to $\sigma$. Let $w_{\sigma,\psi_i}$
be the canonical representative (as defined in \ref{rem:permfacts})
of the element $\rho^{-1}(\psi_i(\sigma)) \in D_n$. Let
\begin{displaymath}
\mathcal{K} = \{[w_{\psi_1,\sigma},w_{\psi_2,\sigma}] | \, \sigma \in A_n \subset S_n\} 
\end{displaymath}
\end{defn}

\begin{lem}[surjectivity of $\rho$ onto $A_n$, when restricted to simple commutators of $D_n$] \label{rho_surj}
The restriction of $\rho$ to the elements of $[D_n,D_n]$ maps surjectively
onto the alternating subgroup $A_n \subset S_n$.
\end{lem}

\noindent \emph{proof}: Given $\sigma \in A_n$, we have that $CR(\sigma) = [\psi_1(\sigma),\psi_2(\sigma)]$
for $\psi_1(\sigma),\psi_2(\sigma) \in S_n$.  This defines an element 

\begin{displaymath}
d  = \rho^{-1}(\psi_1(\sigma))\rho^{-1}(\psi_2(\sigma))\rho^{-1}(\psi_1^{-1}(\sigma))\rho^{-1}(\psi_2^{-1}(\sigma)) \in [D_n,D_n]
\end{displaymath}

Then, since $d \in B_n$ and $\rho$ maps $B_n$ homomorphically to $S_n$, we have that
\begin{displaymath}
\begin{array}{l}
\rho(d) =  \\
 \rho(\rho^{-1}(\psi_1(\sigma))\,\rho^{-1}(\psi_2(\sigma))\,\rho^{-1}(\psi_1^{-1}(\sigma))\,\rho^{-1}(\psi_2^{-1}(\sigma))) = \\
\rho(\rho^{-1}(\psi_1(\sigma)))\,\rho(\rho^{-1}(\psi_2(\sigma)))\,\rho(\rho^{-1}(\psi_1^{-1}(\sigma)))\,\rho(\rho^{-1}(\psi_2^{-1}(\sigma)))  = \\
\psi_1(\sigma) \psi_2(\sigma) \psi_1^{-1}(\sigma) \psi_2^{-1}(\sigma)  = \sigma,
\end{array}
\end{displaymath}
the next to last equality given by the $1-1$-ness of $\rho$ when restricted to $D_n$.\\
\\
We may characterize as follows the set $K$:

\begin{prop}
The set $K \in [D_n,D_n]$ is the set of braids having a braid word 
representative in $\mathcal{K}$. Given an efficient 
algorithm $\mathcal{A}(w)$ for writing any word representative $w$ of $\sigma \in S_n$ 
in canonical form, there is an efficient algorithm to
express any $k \in K$ as a simple commutator of
elements of $D_n$. The algorithm can
accept any $d \in [D_n,D_n]$, and will return $0$ if the
element $d \notin K$.
\end{prop}

\noindent \emph{proof}:  We first show that if $k \in K$, then
it has a braid word representative $w \in \mathcal{K}$. If $k \in K$,
then by definition we have that $CR(\rho(k)) = [\psi_1(\rho(k)),\psi_2(\rho(k))]$
and $[\rho^{-1}(\psi_1(k)),\rho^{-1}(\psi_2(k))] = k$. Then
$w = [w_{\psi_1,\rho(k)},w_{\psi_2,\rho(k)}] \in \mathcal{K}$ and
$w$ is a braid word representative of $k$. \\	
\indent Now if $w = [w_1,w_2] \in \mathcal{K}$, then it represents
a braid  $k \in K$; precisely, the braid $k$ given by $[b_1,b_2]$ where
$b_i$ is the permutation braid canonically represented by $w_i$.  \\
\indent Next, given an element $d \in [D_n, D_n]$, expressed as a braid word
$w_b$ in the standard Artin generators $t_i$ of $B_n$, we can compute its projection
$\rho(d)$, expressed as a word in the standard generators $\tau_i$ of
$S_n$, by replacing each $t_i$ or $t_i^{-1}$ with $\tau_i$. Denote the word so computed by $w_s$. Clearly there is
an algorithm to perform this conversion, such that it will execute on a classical Turing machine, as given above,
in time that is polynomial in the length $l$ of $w_b$ (it will consist of O(l) invocations
of operations of types $1-5$). We can then pass $w_s$ to our implementation of the 
Cejtin-Rivin algorithm, which as we saw is $O(n^2)$ and $O(k)$ where $k$ is the length of $w_s$;
we have that $k$  is bounded by the length of $w_b$. Denote the output of $CR(w_s)$ by $CR(w_s) = [\psi_1,\psi_2]$. The
length of $CR(w_s)$, by definition of our $CR$ implementation, is bounded by a quadratic 
function of $n$. We can compute $\rho^{-1}(\psi_1)$ and $\rho^{-1}(\psi_2)$, expressed
as canonical representatives, by first reducing each  $\psi_i$ to its canonical 
representative $\mathcal{A}(\psi_1)$ (see \ref{subsub:reduce}), and then replacing each $\tau_i$ of 
$\mathcal{A}(\psi_1)$ with $t_i$.  Hence, given an
\emph{arbitrary} braid word representative $w_b$ representing an element of $d \in [D_n,D_n]$,
we can compute the image $w_s = \rho(w_b) \in S_n$, then compute $CR(w_s) = [\psi_1,\psi_2]$,
and then compute $[\rho^{-1}(\mathcal{A}(\psi_1),\rho^{-1}(\mathcal{A}(\psi_2)]$,
all in time  polynomial in the length of $w_b$ (for fixed $n$), and also in $n$ (for fixed input length).

Finally, we can compute $r.g.n.f.(d)$ and $r.g.n.f.([\rho^{-1}(\mathcal{A}(\psi_1),\rho^{-1}(\mathcal{A}(\psi_2)])$. 
The  operation $r.g.n.f.()$, as
noted above, is quadratic in the length of its input. If the
right greedy normal forms of the two words are equal, then $[w_1,w_2]$ expresses
$d$ as a simple commutator, and the output of our algorithm will be the word
$w_1 w_2 w_1^{-1} w_2^{-2}$. If they are not equal, then $d \notin K$, and the
output of our algorithm will be $0$.  Hence our algorithm, as we have just shown, is 
polynomial time in the length of the braid word $w$, representing a braid $b \in B_n$
given as input, and is also polynomial time in $n$.  \indent  $\Box $.

\begin{ex}[Example of computation of $d$ as a simple commutator]

The braid $k = \sigma_2 \sigma_1 \sigma_1 \Omega \in K \subset B_3$ projects
under $\rho$ to the permutation $(123)$. This can be expressed, via
Cejtin-Rivin, as $[\tau_2 \tau_1,\tau_1 \tau_2 \tau_1]$. This
corresponds to $k = [\sigma_2 \sigma_1, \Omega]$.
\end{ex}

\begin{cor} \label{cor: pseudocom}
Let $b \in [B_n, B_n]$ be a braid word representing a product of simple commutators. 
Given an efficient algorithm $\mathcal{A}(w)$ for writing any word representative $w$ of $\sigma \in S_n$ 
in canonical form,  there exists  an efficient algorithm to express $b$
 as the product of a pure braid and a simple commutator of permutation braids.
\end{cor}
\noindent \emph{proof:}  Compute $CR(\rho(b))$, the output of which we denote
$[\psi_1, \psi_2]$. Compute $\mathcal{A}(\psi_1)$ and $\mathcal{A}(\psi_2)$, and
then replace each $\tau_i$ of $\psi_j$ (where $\tau_i$ a generator of $S_n$), with
$t_i$ (the corresponding generator of $B_n$).  This gives an expression $[w_1,w_2]$ for
$d \in [D_n,D_n]$, in a canonical form induced by $\mathcal{A}$,  by the surjectivity of 
$\rho:[D_n,D_n] \rightarrow A_n \subset S_n$.  Since the kernel of $\rho$ is the subgroup of pure
braids $P \subset B_n$, we have that $b = p\,d$ for some $p \in P$, where
$p = b[d_1,d_2]^{-1}$.  By the preceding proposition and the conditions of
our hypothesis,  this expression for $b$ can be computed in time polynomial
in the length of $b$.  Hence there is an efficient algorithm to express
any element of the commutator subgroup of $B_n$ as the product of a pure braid and a simple
commutator of permutation braids. \indent $\Box$

%\vspace*{3.0in}
\pagebreak

\subsection{A probablistic algorithm for the SCPP for $B_n$}

Now, given a braid word $b$ representing an element of $\gamma \in B_n$, such that we are guaranteed
$\gamma$ is a simple commutator, we can write $b$ in the form $b = p[d_1,d_2]$ as above.  We know
there must be a sequence $p[d_1,d_2] = w_0 \rightarrow w_1 \rightarrow  \cdots \rightarrow w_n = [a,b]$
where each $w_i$ is obtained from $w_{i-1}$ by one subword substitution corresponding to a defining
relation of the braid group. We know moreover that there is some minimal $N_b$ such that there is a sequence
$p[d_1,d_2] = v_0 \rightarrow v_1 \rightarrow  \cdots \rightarrow v_n = [a,b]$, having
a subset of terms $v_0 = v_{k_1}, v_{k_2}, \cdots ,v_{k_m} = v_n$,
where $1 \leq k_1 < k_2 < \cdots <k_m \leq n$,  such that for all $i$:\\
\\
\indent \indent (1) $ k_{i+1} - k_i < N_b$\\
\indent \indent (2) $v_i = h_i[a_i, b_i]$, where $| h_i | < | p_i]$, for $h_i, a_i, b_i \in B$. \\
\indent \indent (3) $|h_i+1| < |h_{i}|$ \\
\indent  \indent (4) $|h_1| < |p|$ \\
\indent \indent  for $h_i, a_i, b_i \in B$, where $| |$ denotes word length. \\

We accordingly define the $K$-bounded version of the SCPP problem as
follows: given $b$ an arbitrary braid word  representing a simple
commutator, find such a sequence, for $N_b = K$, if one exists. We
can define a probablistic algorithm for
searching for such a sequence as follows:\\
\\
INPUT: A braid word $b$ guaranteed to be a simple commutator. \\
\\
ALGORITHM: \\
\\
1. Randomly select integer $M \in \mathbf{Z^+}$.  Write the word $b$ in the form $b_0 = p[d_1,d_2]$ in the
manner described above.  Let the output $O$ be set to $e$ (the empty word). \\
\noindent 2. For each $i = 0,  1, \cdots M$: \\
\indent 3. Randomly select $N_i \in \mathbf{Z^+}$ s.t. if $i > 0$, then $N_i > N_{i-1}$. Randomly Let $j = 0$.  \\
\indent 4.   While $j < N_i$:\\
\indent \indent 5.  Randomly select a  subword of $b_i$ that can be rewritten via a rule given by one of the
braid relations, and apply the rule. \\
\indent \indent 6. If the result $b_{i+1}$ is of the form $[a,b]$, append  the symbol $\rightarrow$ followed
by the word $[a,b]$ followed by ($N_i$) to $O$. Terminate with success. \\
\indent \indent 7. Else if $b_i$ satisfies properties $(2)-(4)$ above, append the symbol $\rightarrow$ followed
by the word $b_{i+1}$ to $O$. set $j=0$. \\
\indent \indent 8. Else $j = j+1$. \\
9. Terminate with failure. \\

\pagebreak

\section{Summary and questions for further research}

In this manuscript, we first discussed motivations for studying the
commutator subgroups of various finite groups. We then defined the
simple commutator promise problem (SCPP), and studied the problem
for the case of the symmetric group $S_n$, applying the work of
Cejtin and Rivin. We then extended the algorithm to a small but
noteable subset of elements of the braid group $B_n$; namely, to a
subset $K$ of the set of simple commutators of permutation braids,
such that $K$ is in $1-1$ correspondence with the simple commutators
of $S_n$. We found that any element of the commutator subgroup of $B_n$ can be
expressed efficiently as the product of an element of $K$ and a pure braid. \\
\\
We here pose a few questions for further research: \\

\noindent $1$.  Can the Cejtin-Rivin algorithm be extended to a larger
subset of braids than the set $K$ defined above (\ref{def:kwords})? What
are the necessary and sufficient conditions that the commutator subgroup $[G,G]$ of 
a given group $G$ must satisfy, in order for the simple commutator promise problem to be 
efficiently solveable for $G$? \\
\\
$2$. Given an arbitrary braid word $b$ that represents a simple commutator of braids, are there
values of $N_b$ and $M_b$ such that the probabilistic algorithm of the preceeding section will terminate with a 
probability of success $P$, such that $P > \frac{1}{2}$, and such that (for fixed $n$) the algorithm will execute 
on a single-tape classical deterministic Turing machine in time polynomial in the length of the input $b$, and such
that the algorithm will execute in time polynomial in $n$, for input of fixed length? \\
\\
$3$. The above discussion of complexity concerns worst-case scenarios, but there are other measures of
complexity; see for example the description of average-case complexity in \cite{Kapovich}. Can we use the a
above results to find an algorithm that solves the SCDP for $B_n$, with average-case complexity polynomial in both $n$ and the length of the algorithm input, relative to a discrete probability measure $\mu$ on the set of braid words? What measure
should be used? \\
\\
$4$. It is known that braid groups have a faithful linear
representation \cite{Bigelow} \cite{Budney}. Can this representation be used to
obtain a solution to the simple commutator promise problem for $B_n$? \\
\\
$5$. There are alternatives to the classical Turing machine model of
computing, in particular the \emph{quantum} model, which replaces
the transition functions of the Turing machine with unitary
operations on a finite dimensional Hilbert space, as described in
\cite{Nielsen}. It is known that braid groups can be used as a model
for a method of computation that is equivalent in power to the
quantum model
\cite{FreedmanA}\cite{FreedmanB}\cite{FreedmanB}\cite{FreedmanD}\cite{FreedmanE}.
Can either of these models yield an efficient solution to the simple commutator promise problem 
for the braid groups?

\pagebreak

\section{Figures}

\begin{picture}(150,250)
\thicklines

\put(40,150) {\line(5,1) {100}}
\put(140,171){\line(1,0) {200}}
\put(40,150) {\line(1,0) {200}}
\put(241,150) {\line(5,1){100}}
\put(40,1) {\line(5,1) {100}}
\put(140,21) {\line(1,0) {200}}
\put(40,1) {\line(1,0) {200}}
\put(241,1){\line(5,1) {100}}
\put(100,155){\line(1,-6){25}}
\put(130,155){\line(1,-6){25}}
\put(160,155){\line(1,-6){12}}
\put(174,70){\line(1,-6){11}}
\put(190,155){\line(-1,-6){12}}
\put(178,82){\line(-3,-4){20}}
\put(142,42){\line(-1,-1){15}}
\put(116,20){\line(-1,-2){8}}

\put(230,70){\vector(0,1){40}}
\put(230,70){\vector(-1,-2){15}}
\put(230,70){\vector(1,0){40}}
\put(230,112){z}
\put(272,70){y}

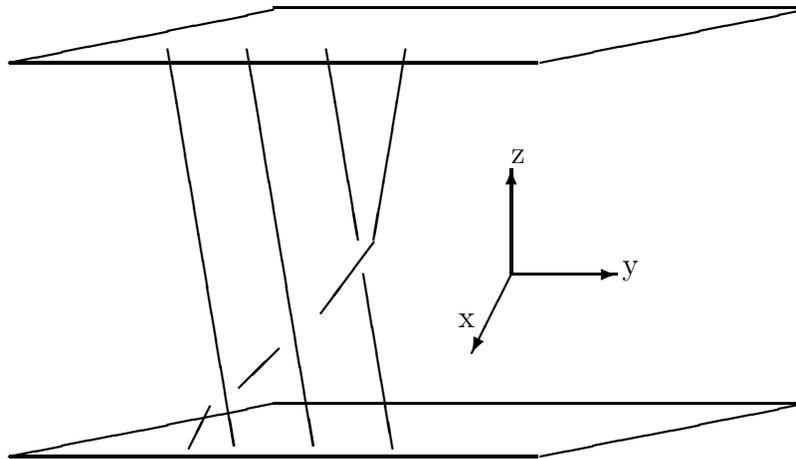
\put(1,-14){figure 1: a braid diagram corresponding (bottom to top) to $\sigma_1 \sigma_2 \sigma_3^{-1}$}
\put(210,50){x}

\end{picture}

\begin{picture}(150,250)
\thicklines

\put(50,1) {\line(0,1){100}}
\put(90,1) {\line(0,1){100}}
\put(110,50) {\circle*{2}}
\put(130,50) {\circle*{2}}
\put(150,50) {\circle*{2}}
\put(170,1){\line(0,1){30}}
\put(200,1){\line(0,1){30}}
\put(170,31){\line(1,1){10}}
\put(200,31){\line(-1,1){30}}
\put(190,48){\line(1,1){12}}
\put(170,61){\line(0,1){39}}
\put(200,61){\line(0,1){39}}
\put(220,50) {\circle*{2}}
\put(240,50) {\circle*{2}}
\put(260,50) {\circle*{2}}
\put(280,1) {\line(0,1){100}}
\put(48,-8){1}
\put(88,-8){2}
\put(168,-8){i}
\put(198,-8){i+1}
\put(278,-8){n}
\put(48,-20){figure 2: a braid diagram corresponding to $\sigma_i$}
\end{picture}

\begin{picture}(150,250)
\put(109,55){\line(1,-1){12}}
\put(109,25){\line(1,1){30}}
\put(130,37){\line(1,-1){20}}
\put(109,1){\line(1,0){30}}
\put(140,1){\line(1,1){54}}
\put(156,10){\line(1,-1){10}}
\put(140,55){\line(1,0){31}}
\put(172,55){\line(1,-1){10}}
\put(185,42){\line(1,-1){10}}
\put(166,0){\line(1,0){35}}
\put(40,-12){figure 3: a braid diagram corresponding to $\Omega \in B_3$}
\end{picture}
\pagebreak

%\begin{figure}[h]
% \includegraphics[3in,8in][6in,10in]{figure1.jpg} \caption{a braid diagram corresponding (bottom to top) to $\sigma_1 %\sigma_2 \sigma_3^{-1}$, and its/
%  inverse}
% \label{fig:1}
%\end{figure}
%
%\begin{figure}[h]
% \includegraphics[3in,6in][6in,8.5in]{figure2.jpg} \caption{a braid diagram corresponding to $\sigma_i$}
% \label{fig:2}
%\end{figure}
%
%\pagebreak
%
%\begin{figure}[h]
% \includegraphics[0in,0in][2in,2in]{figure3.jpg} \caption{a braid diagram corresponding to $\Omega \in B_3$}
% \label{fig:3}
%\end{figure}

%\pagebreak
% ------------------------------------------------------------------------
%GATHER{Xbib.bib}   % For Gather Purpose Only
%GATHER{Paper.bbl}  % For Gather Purpose Only
\bibliographystyle{amsplain}

\normalsize

\normalfont

\linespread{1}
\pagebreak

%\section{Acknowledgements}
%I thank Professors Michael Anshel and Brendan Owens for their
%assistance in the preparation of this proposal. Of course, all
%errors it contains  are my own.
\end{document}